\documentclass[a4paper,10pt]{article}
\usepackage{amscd}
\usepackage{bbding}
\usepackage{amsfonts}
\usepackage{graphicx,amssymb,amsthm}
\usepackage{latexsym,amsmath,amscd,cite,enumerate}
\usepackage{float}
\usepackage{leftidx}
\usepackage{epstopdf}
\usepackage{manfnt}
\usepackage{stmaryrd}
\usepackage{cases}
\usepackage{bbm}
\usepackage{tabls}
\usepackage{multirow}
\usepackage{color}
\usepackage{tensor}

\allowdisplaybreaks[4]
\DeclareFontEncoding{OT2}{}{}
\DeclareTextSymbol{\cyrsftsn}{OT2}{126}
\DeclareTextSymbol{\textnumero}{OT2}{125}

\theoremstyle{definition}
\newtheorem{theorem}{Theorem}[section]
\newtheorem{lemma}{Lemma}[section]
\newtheorem{corollary}{Corollary}[section]

\newtheorem{definition}{Definition}[section]
\newtheorem{remark}{Remark}[section]

\begin{document}

\title{{A new nonlocal fractional differential quasi-variational inequality in Hilbert spaces with applications}\thanks{This work was supported by the National Natural Science Foundation of China (11901273, 62272208, 12171339, 12171070), the Program for Science and Technology Innovation Talents in Universities of Henan Province (23HASTIT031), the Young Backbone Teachers of Henan Province (2021GGJS130), the Natural Science foundation of Sichuan Province (2024NSFSC1392).}}
\author{Zeng-bao Wu$^{a}$, Tao Chen$^{b}$, Quan-guo Zhang$^{a}$, Yue Zeng$^{c}$, Nan-jing Huang$^{c}$\footnote{Corresponding author. E-mail addresses: njhuang@scu.edu.cn, nanjinghuang@hotmail.com}, Yi-bin Xiao$^{d}$\\
$^{a}${\small \textit{Department of Mathematics, Luoyang Normal University, Luoyang, Henan  471934, P.R. China}}\\
$^{b}${\small \textit{School of Science, Southwest Petroleum University, Chengdu, Sichuan 610500, P.R. China}}\\
$^{c}${\small \textit{Department of Mathematics, Sichuan University, Chengdu, Sichuan 610064, P.R. China}}\\
$^{d}${\small \textit{School of Mathematical Sciences, University of Electronic Science and Technology  of China, }}\\
{\small\textit{Chengdu, Sichuan 611731, P.R. China}}
}
\date{ }
\maketitle

\begin{flushleft}
\hrulefill\newline
\end{flushleft}
\textbf{Abstract}.  This paper considers a new nonlocal fractional differential quasi-variational inequality (NFDQVI) comprising a fractional differential equation with a nonlocal condition and a time-dependent quasi-variational inequality in Hilbert spaces. Qualitative properties of the solution for the time-dependent parameterized  quasi-variational inequality are investigated, which improve some known results in the literature. Moreover, the unique existence of the solution and Hyers-Ulam stability are obtained for such a novel NFDQVI under mild conditions.  Finally,  the obtained abstract results for NFDQVI are applied to analyze the unique solvability and stability addressing a time-dependent multi-agent optimization problem and a time-dependent price control problem.
\begin{flushleft}
\hrulefill
\end{flushleft}

\section{Introduction}
\noindent \setcounter{equation}{0}
As is well known, differential variational inequalities (DVIs for brevity), investigated in Euclidean spaces by Pang and Stewart \cite{PS} for the first time, are a class of coupled systems comprising variational inequalities (VIs for short) and differential equations which generalize the notion of differential algebraic equations, evolution VIs, projected dynamical systems, differential complementarity problems, etc. Due to the widespread applications in the fields of dynamic transportation network, ideal diode circuits, microbial fermentation processes, contact problems, dynamic Nash equilibrium problems, price control problems and so on, DVIs in finite/infinite dimensional spaces have attracted more and more researchers to discuss the theory, algorithms and applications. For instance, Wang et al. \cite{WQTX} studied the solvability of a delay DVI, and they also performed the convergence analysis on the provided algorithm. Zhang et al. \cite{zhang2023s} investigated a numerical approximation method for solving a stochastic DVI and provided some applications in stochastic environment. Mig\'{o}rski et al. \cite{MCD} obtained the solvability of a differential variational-hemivariational inequality and  provided an application to contact mechanics. For more works, we refer the readers to the excellent review on DVIs by Brogliato and Tanwani \cite{BT} and the citations therein.

It is worth pointing out that fractional calculus (FC for brevity) provide a better tool than integer order ones to describe many phenomena and physical processes, especially those with hereditary and memory properties, which has been widely used in various fields such as biology, mechanics and dynamic systems, signal and image processing, environmental science, economics, materials and so on (see, e.g., \cite{SZBCC,NFAC,ZSSZH,YWZ}). In 2015, Ke et al. \cite{K2015} combined FC with DVIs for the first time, and they investigated a fractional DVIs with delay. In 2021, Wu et al. \cite{WWHWW} studied a fuzzy fractional DVI consisting of a VI and a fuzzy fractional differential inclusion, and they \cite{WLZX} in 2024 further investigated the stability of solutions for the fuzzy fractional DVI. Recently, Zeng et al. \cite{ZZH} studied the unique solvability of a stochastic fractional DVI and applied the result to a spatial price equilibrium problem in stochastic environments. For more related results, the readers are encouraged to consult \cite{LKOZ,WCH,WCLH,WWHXZ,ZBYN,ZCN,WC2024} and the citations therein.

Since nonlocal conditions are more effective than classical initial conditions in physical applications, nonlocal Cauchy problems (particularly fractional-order ones) have received growing research interest (see, e.g., \cite{AL,LB1991,BP2007, BP2009,DW2016,DWZ,WIF} and the reference therein).  Recently, Feng et al. \cite{FYZ} studied a generalized time-fractional super-diffusion equation with a nonlocal integral observation. Eckardt and Zhigun \cite{EZ2025} investigated the solvability of a nonlocal equation with degenerate anisotropic diffusion. Chen et al. \cite{CZHX} examined a new nonlocal impulsive fractional differential hemivariational inclusions and provided an application to a frictional contact problem. We would like to mention that differential quasi-variational inequalities (DQVIs for brevity), as one of the important research contents of DVIs, have attracted the attention of more and more researchers due to their wide application in the generalized differential Nash game, frictional contact problems and so on (see, e.g., \cite{DLJX,PS,WTLH,WCH}). Recently, Chu et al. \cite{CCHX} studied the unique solvability of a novel DQVI and applied their findings to a frictional contact problem. Zhao et al. \cite{ZCL} examined the existence of solutions for a differential quasi-variational-hemivariational inequality. However, to the best of our knowledge, there are few works addressing the nonlocal fractional DQVI (NFDQVI for short). Thus, it would be interesting and important both in theory and in practice to investigate NFDQVIs under some mild conditions.

On the other hand, Hyers-Ulam (H-U for brevity) stability refers to a theorem's argument being true or almost true if we make minor changes to its assumptions, which is a fundamental type of stability for functional equations. This type of concept has been generalized to various classes of ordinary and partial differential equations, (set-valued) functional equations, and the study of this region has become one of the most significant issues in the field of mathematical analysis (see, e.g., \cite{AXZ,Jung,Jung15,POS}). In particular, in 2011, Wang et al. \cite{WLZ} first introduced the notion of H-U stability into fractional differential equations (FDEs for brevity). In 2012, Wang et al. \cite{WLZ12,WZ} continued to study the H-U stability of FDEs by using a fixed point approach, and also investigated the Mittag-Leffler-Ulam stability of fractional evolution equations.  Recently, some concepts of H-U stability have already been introduced to investigate the stability of DVIs. For example,  Loi and Vu \cite{LV} first studied the H-U stability and the unique existence of solution for a DVI with nonlocal conditions and provided two applications. Du et al. \cite{DLJX} established the H-U stability and unique solvability for a fractional DQVI with Katugampola fractional operator. Jiang et al. \cite{JSZ} examined the H-U stability and unique solvability for a random DVIs with nonlocal boundary conditions and provided two applications. However, the study of H-U stability for NFDQVIs is still not fully explored and much is desired to be done.

The present paper is thus devoted to the study of the unique solvability and the H-U stability for the following NFDQVI
\begin{equation}\label{GFDQVI}
\left\{
\begin{array}{l}
\tensor*[^{C}_0]{\mathrm{D}}{_s^q} x(s) = f(s,x(s)) + g(s,x(s),u(s)), \; \forall  \, s \in I,\\
u(s) \in \mbox{SOL}\left(K(\cdot),G(s,x(s),\cdot)\right), \; \forall  \,  s \in I, \\
x\left(0\right)=x_{0}+\psi(x),
\end{array}
\right.
\end{equation}
where $\tensor*[^{C}_0]{\mathrm{D}}{_s^q}$ is the Caputo type fractional derivative with $q \in (0,1]$, $f:I \times \mathcal{H}_1 \rightarrow \mathcal{H}_1$ and $g:I \times \mathcal{H}_1 \times \mathcal{H}_2 \rightarrow \mathcal{H}_1$ are two given mappings with $I=[0,T]$, $\mathcal{H}_1$ and $\mathcal{H}_2$ being two Hilbert spaces, $K: \Omega \rightarrow 2^{\Omega}$ is a set-valued mapping with nonempty convex and closed values, $G:I \times \mathcal{H}_1 \times \Omega \rightarrow \mathcal{H}_2$ is a given mapping, $\Omega (\neq \emptyset) \subset \mathcal{H}_2$ is a convex and closed set, $x_{0} \in \mathcal{H}_1$, $\psi: C(I,\mathcal{H}_1)  \rightarrow \mathcal{H}_1$ is a given mapping and $\mbox{SOL}\left(K(\cdot),G(s,z,\cdot)\right)$ is the solution set to the following time-dependent parameterized quasi-variational inequality (PQVI for brevity): for any given  $s \in I$ and $x(s)\in \mathcal{H}_1$, find $u(s) \in K(u(s))$ such that
\begin{equation}\label{PQVI}
  \left<G(s,x(s),u(s)),v-u(s)\right> \geq 0, \; \forall \; v \in K(u(s)),
\end{equation}
where $\left<\cdot,\cdot\right>$ is the inner product in $\mathcal{H}_2$. Clearly, if $K(u)=\phi(u)+\Omega$ with $\phi: \Omega \rightarrow \Omega$ and $\Omega$ being a convex and closed cone of $\mathcal{H}_2$, then PQVI (\ref{PQVI}) is equivalent to the time-dependent parameterized quasi-complementarity problem (PQCP for brevity): for any given  $s \in I$ and $x(s)\in \mathcal{H}_1$,  find $u(s) \in \phi(u(s))+\Omega$ such that
\begin{equation}\label{PQCP}
 G(s,x(s),u(s)) \in  \Omega^*, \quad   \left<G(s,x(s),u(s)),u(s)-\phi(u(s))\right> =0,
\end{equation}
where $\Omega^*=\{\omega  \in \mathcal{H}_1: \left<\omega ,v\right>\geq 0, \; \forall v \in \Omega\}$. It is worth mentioning that NFDQVI (\ref{GFDQVI}) includes many kinds of DVIs as its special cases. For example:
\begin{itemize}
\item[(i)] If $f=0$, then (\ref{GFDQVI}) can be rewritten as
\begin{equation*}
\left\{
\begin{array}{l}
\tensor*[^{C}_0]{\mathrm{D}}{_s^q} x(s) = g(s,x(s),u(s)), \; \forall  \, s \in I,\\
u(s) \in \mbox{SOL}\left(K(\cdot),G(s,x(s),\cdot)\right), \; \forall  \,  s \in I, \\
x\left(0\right)=x_{0}+\psi(x),
\end{array}
\right.
\end{equation*}
which is an emerging problem.
  \item[(ii)] If $K(u)=\phi(u)+\Omega$ with $\phi: \Omega \rightarrow \Omega$ and $\Omega$ being a convex and closed cone, then (\ref{GFDQVI}) is becoming to a nonlocal fractional differential quasi-complementarity problem (NFDQCP for short) as follows:
\begin{equation}\label{GFDQCP}
\left\{
\begin{array}{l}
\tensor*[^{C}_0]{\mathrm{D}}{_s^q} x(s) = f(s,x(s)) + g(s,x(s),u(s)), \; \forall  \, s \in I,\\
G(s,x(s),u(s)) \in  \Omega^*, \;   \left<G(s,x(s),u(s)),u(s)-\phi(u(s))\right> =0, \; \forall  \,  s \in I, \\
x\left(0\right)=x_{0}+\psi(x),
\end{array}
\right.
\end{equation}
which is also an emerging problem.
  \item[(iii)] If $q=1$, $\mathcal{H}_1=R^n$, $\Omega=\mathcal{H}_2=R^m$, $\psi(x)=\mathbf{0}_{R^n}$ with $\mathbf{0}_{R^n}$ being the zero vector in $R^n$, $g(s,x(s),u(s))=g_1(s,x(s))\cdot u(s)$, $g_1:I \times R^n \rightarrow R^{n \times m}$, $G(s,x(s),u(s))=G_1(s,x(s))+ G_2(u(s))$, $G_1:I \times R^n \rightarrow R^m$, and $G_2: \Omega\rightarrow R^m$, then (\ref{GFDQVI}) is becoming to a specific case of the DQVI, which was considered by Wang et al. \cite{WTLH}.
  \item[(iv)] If $q=1$, $\mathcal{H}_1=R^n$, $\mathcal{H}_2=R^m$, $\psi(x)=\mathbf{0}_{R^n}$, $K(x)=\Omega$ is independent of $x$, $g(s,x(s),u(s))=g_1(s,x(s))\cdot u(s)$, $g_1:I \times R^n \rightarrow R^{n \times m}$, $G(s,x(s),u(s))=G_1(s,x(s))+ G_2(u(s))$, $G_1:I \times R^n \rightarrow R^m$, and $G_2: \Omega\rightarrow R^m$,  then (\ref{GFDQVI}) is reduced to the DVI, which was systematically investigated by Pang and Stewart \cite{PS}.
\end{itemize}

The main contributions of the current work are three-fold. First, interesting properties of the solution for PQVI in (\ref{GFDQVI}) are investigated under the hypotheses of strong pseudomonotonicity and Lipschitzean and some sufficient conditions are also given for ensuring the unique solvability of NFDQVI (\ref{GFDQVI}) by using the properties of the solution for PQVI and the Banach fixed point theorem. Second, two H-U stability results for NFDQVI (\ref{GFDQVI}) are obtained under some mild conditions by carefully calculations and notably employing the inequalities techniques. Last but not least, applications of the abstract results are given to a time-dependent multi-agent optimization problem and time-dependent a price control problem, respectively.

The rest of this paper is organized as follows. Next section recalls some notations, definitions, and lemmas. In Section 3, some properties of solution for PQVI in (\ref{GFDQVI}) are examined under the hypotheses of strong pseudomonotonicity and Lipschitzean. Moreover, the unique solvability for NFDQVI (\ref{GFDQVI}) is proved based on the Banach fixed point principle. In Section 4, two H-U stability results for NFDQVI (\ref{GFDQVI}) are provided by tedious calculations and notably utilising inequalities techniques.  As applications, the unique solvability and stability for a time-dependent multi-agent optimization problem and time-dependent a price control problem are investigated in Section 5. At last, the conclusion is given in Section 6.

\section{Preliminaries}
\noindent \setcounter{equation}{0}
In this section, we review some useful knowledge in FC and variational analysis. Let $R_+$ be the set of nonnegative reals. Let $Y$ be a Banach space with the norm $\|\cdot\|$ and $2^Y$ stand for the set of all subsets of $Y$. Let $C(I,Y)=\{f:I\to Y: f \mbox{ is  continuous }\}$ be a Banach space with Bielecki's norm
$$\|x\|_\mathcal{B}=\underset{s \in I} \max \; e^{-\gamma s}\|x(s)\|, \;\; x \in C(I,Y).$$
It is easy to see that the above norm is equivalent to usual norm in $C(I,Y)$ (see, e.g., \cite[Theorem 8]{ZCN}), where $\gamma>0$ is a given constant. From now on, $P_{\Omega_0}$ denotes the projection of a real Hilbert space $\mathcal{H}$ onto the set $\Omega_0 \subset \mathcal{H}$.

\begin{definition}\cite{KD}
Let $0<q \leq 1$.
\begin{itemize}
  \item[(i)] The Riemann-Liouville fractional integral of order $q$ for the function $x$ is given by
  \begin{equation*}
\tensor*{\mathrm{I}}{_0^q}x(s) = \frac{1}{\Gamma(q)}\int_0^s(s-\varsigma)^{q-1}x(\varsigma)d\varsigma,
\end{equation*}
where the Gamma function $\Gamma$ is defined by $\displaystyle \Gamma(q)=\int_{0}^\infty \varsigma^{q-1}e^{-\varsigma}d\varsigma$.
  \item[(ii)] The Riemann-Liouville fractional derivative of order $q$ for the function $x$ is given by
      \begin{eqnarray*}
\tensor*[_0]{\mathrm{D}}{_s^q}x(s) =
\frac{1}{\Gamma(1-q)}\frac{d}{ds}\int_0^s(s-\varsigma)^{-q}x(\varsigma)d \varsigma.
\end{eqnarray*}
  \item[(iii)] The Caputo fractional derivative of order $q$ for the function $x$ is defined by
      \begin{eqnarray*}
\tensor*[^{C}_0]{\mathrm{D}}{_s^q}x(s) =
\tensor*[_0]{\mathrm{D}}{_s^q}\left(x(s)-x(0)\right).
\end{eqnarray*}
\end{itemize}
\end{definition}

\begin{remark}
\begin{itemize}
   \item[(i)] If $q=1$, then $\tensor*[^{C}_0]{\mathrm{D}}{_s^q}x(s)$ is classical derivative $x^\prime(s)$.
   \item[(ii)] If $x$ is a continuous differentiable function on $I$, then
\begin{eqnarray*}
\tensor*[^{C}_0]{\mathrm{D}}{_s^q}x(s) =
\frac{1}{\Gamma(1-q)}\int_0^s(s-\varsigma)^{-q}x^\prime(\varsigma)d \varsigma.
\end{eqnarray*}
  \item[(iii)] The integrals that show up in the definitions above are taken in Bochner's sense if $f$ is an abstract function with values in Banach space $X$.
\end{itemize}
\end{remark}

\begin{definition}
Let $\Omega_0\neq\emptyset$ be a set of a real Hilbert space $\mathcal{H}$. We call that $G: \Omega_0 \rightarrow \mathcal{H}$ is
\begin{itemize}
  \item[(a)] strongly monotone on $\Omega_0$ if there is $\eta>0$ such that
   $$\left<G(w_2)-G(w_1),w_2-w_1\right> \geq \eta\|w_2-w_1\|^2 \quad \mbox{for all} \;\; w_1, w_2 \in \Omega_0;$$
  \item[(b)] strongly pseudomonotone on $\Omega_0$ if there is $\eta>0$ such that
   $$\left<G(w_2),w_1-w_2\right> \geq 0 \;\; \Rightarrow \;\; \left<G(w_1),w_1-w_2\right> \geq \eta\|w_1-w_2\|^2 \quad \mbox{for all} \;\; w_1, w_2 \in \Omega_0.$$
\end{itemize}
\end{definition}

Obviously, the implication holds: (a) implies (b). In general, the converse implications are not true.

Motivated by \cite[Theorem 4.1]{BP2009} and \cite[Lemma 2.3]{WCH}, we introduce the definition of NFDQVI (\ref{GFDQVI}) as follows.

\begin{definition}\label{def-FDQVI} For $(x,u)\in C(I,\mathcal{H}_1)\times C(I,\Omega)$, we call that the pair $(x,u)$ is a solution of NFDQVI (\ref{GFDQVI}) if
\begin{equation} \label{eq-fg}
\left\{
\begin{array}{l}
x(s) = x_{0} + \psi(x)+\frac{1}{\Gamma(q)} \int_0^s(s-\zeta)^{q-1}\left[f(\zeta,x(\zeta)) + g(\zeta,x(\zeta),u(\zeta))\right]d\zeta, \; \forall  \, s \in I,\\
u(s) \in \mbox{SOL}\left(K(\cdot),G(s,x(s),\cdot)\right), \; \forall  \,  s \in I.
\end{array}
\right.
\end{equation}
Within it, we say that $x$ is the trajectory and $u$ is the variational control trajectory.
\end{definition}

Motivated by \cite[Definition 2]{LV} and \cite[Definitions 3.1-3.4]{WZ}, we introduce the Mittag-Leffler-Hyers-Ulam (MLHU for short) stability  concepts for NFDQVI (\ref{GFDQVI}) as follows.

\begin{definition}\label{def-MLUH}
NFDQVI (\ref{GFDQVI}) is called MLHU stable w.r.t. the function $\mathrm{E}_q$, if there exists $c>0$ such that for each number $\varepsilon>0$ and solution $(z,\upsilon)\in C(I,\mathcal{H}_1)\times C(I,\Omega)$ of the following inequality system
\begin{equation}
\left\{
\begin{array}{l}
\left\|\tensor*[^{C}_0]{\mathrm{D}}{_s^q} z(s) - f(s,z(s)) - g(s,z(s),\upsilon(s))\right\| \leq \varepsilon, \; \forall  \, s \in I,\\
\upsilon(s) \in \mbox{SOL}\left(K(\cdot),G(s,z(s),\cdot)\right), \; \forall  \, s \in I, \\
z\left(0\right)=x\left(0\right),
\end{array}
\right.  \label{MLUH}
\end{equation}
there exists a solution $(x,u)$ of NFDQVI (\ref{GFDQVI})  it holds
$$\|z(s)-x(s)\| \leq  c \varepsilon \mathrm{E}_q(s), \; \forall \, s \in I.$$
\end{definition}

\begin{definition}\label{def-GMLUH}
NFDQVI (\ref{GFDQVI}) is called generalized MLHU stable w.r.t. the function $\mathrm{E}_q$, if there exists $\theta \in C(R_+,R_+)$ with $\theta(0)=0$, such that for each solution $(z,\upsilon)\in C(I,\mathcal{H}_1)\times \in C(I,\Omega)$ of (\ref{MLUH}) there exists a solution $(x,u)$ of NFDQVI (\ref{GFDQVI}) with
$$\|z(s)-x(s)\| \leq  \theta(\varepsilon)  \mathrm{E}_q(s), \; \forall \,  s \in I.$$
\end{definition}

\begin{definition}\label{def-MLUHR}
For any given $\varphi \in C(I,R_+)$, NFDQVI (\ref{GFDQVI}) is called Mittag-Leffler-Hyers-Ulam-Rassias (MLHUR for short) stable w.r.t. the function $\varphi\mathrm{E}_q$, if there exists $c_\varphi>0$ such that for each $\varepsilon>0$ and each solution $(z,\upsilon) \in C(I,\mathcal{H}_1)\times  C(I,\Omega)$ of the following inequality system
\begin{equation}
\left\{
\begin{array}{l}
\left\|\tensor*[^{C}_0]{\mathrm{D}}{_s^q} z(s) - f(s,z(s)) - g(s,z(s),\upsilon(s))\right\| \leq \varepsilon \varphi(s), \; \forall  \, s \in I,\\
\upsilon(s) \in \mbox{SOL}\left(K(\cdot),G(s,z(s),\cdot)\right), \; \forall  \, s \in I, \\
z\left(0\right)=x\left(0\right),
\end{array}
\right.   \label{MLUHR}
\end{equation}
there exists a solution $(x,u)$ of NFDQVI (\ref{GFDQVI}) with
$$\|z(s)-x(s)\| \leq  c_\varphi \varepsilon \varphi(s) \mathrm{E}_q(s), \; \forall \,  s \in I.$$
\end{definition}

\begin{definition}\label{def-GMLUHR}
For any given $\varphi \in C(I,R_+)$, NFDQVI (\ref{GFDQVI}) is called generalized MLHUR stable w.r.t. the function $\varphi \mathrm{E}_q$, if there exists a real number $c_\varphi>0$ such that for each solution $(z,\upsilon) \in C(I,\mathcal{H}_1)\times  C(I,\Omega)$  of the following inequality system
\begin{equation}
\left\{
\begin{array}{l}
\left\|\tensor*[^{C}_0]{\mathrm{D}}{_s^q} z(s) - f(s,z(s)) - g(s,z(s),\upsilon(s))\right\| \leq  \varphi(s), \; \forall  \, s \in I,\\
\upsilon(s) \in \mbox{SOL}\left(K(\cdot),G(s,z(s),\cdot)\right), \; \forall  \, s \in I, \\
z\left(0\right)=x\left(0\right),
\end{array}
\right.  \label{GMLUHR}
\end{equation}
there exists a solution $(x,u)$ of NFDQVI (\ref{GFDQVI}) it holds
$$\|z(s)-x(s)\| \leq  c_\varphi \varphi(s) \mathrm{E}_q(s), \; \forall \,  s \in I.$$
\end{definition}

\begin{remark}\label{ISys}
For $(z,\upsilon) \in C(I,\mathcal{H}_1)\times C(I,\Omega)$, the pair $(z,\upsilon)$ is a solution of the inequality system (\ref{MLUH}) if and only if there exists a function $h\in C(I, \mathcal{H}_1)$ such that
\begin{itemize}
  \item[(i)] for any $s\in I$, $\|h(s)\|\leq \varepsilon$;
  \item[(ii)] $(z,\upsilon)$ is the solution of FDQVI 
$$
\left\{
\begin{array}{l}
\tensor*[^{C}_0]{\mathrm{D}}{_s^q} z(s) = f(s,z(s)) + g(s,z(s),\upsilon(s)) + h(s), \; \forall  \, s \in I,\\
\upsilon(s) \in \mbox{SOL}\left(K(\cdot),G(s,z(s),\cdot)\right), \; \forall  \,  s \in I, \\
z\left(0\right)=x\left(0\right).
\end{array}
\right.
$$
\end{itemize}
\end{remark}

We have similar remarks for the inequality systems (\ref{MLUHR}) and (\ref{GMLUHR}).

\begin{remark}
According to Definitions \ref{def-MLUH}-\ref{def-GMLUHR}, the following implications hold: (i) Definition \ref{def-MLUH} implies Definition \ref{def-GMLUH}; (ii) Definition \ref{def-MLUHR} implies  Definition \ref{def-GMLUHR}; (iii) Definition \ref{def-MLUHR} implies Definition \ref{def-MLUH}.
\end{remark}

\begin{lemma}\label{POmega}\cite[Theorem 3.16]{BC}  Let $\Omega_0\neq\emptyset$ be a convex and closed set of a real Hilbert space $\mathcal{H}$ endowed with inner product $\left<\cdot,\cdot\right>$. Then it holds
$$w=P_{\Omega_0}\left[\varpi\right] \;\; \Leftrightarrow \;\; w \in \Omega_0 \;\; \mbox{and}\;\; \left<w-\varpi,v-w\right> \geq 0, \; \forall \; v \in \Omega_0,$$
where $w,\varpi \in \mathcal{H}$.
\end{lemma}

\begin{lemma}\cite[Theorem 5.1]{CP82}\label{QVI}
Given a set-valued mapping $K:\mathcal{H} \rightarrow 2^\mathcal{H}$ with convex and closed values, and a nonlinear operator $F:\mathcal{H} \rightarrow \mathcal{H}$ with $\mathcal{H}$ being a real Hilbert space. Then $u \in \mathcal{H}$ solves quasi-variational inequality $\mbox{QVI}(K(\cdot), F)$: find $w \in K(w)$ such that $\left<F(w),v-w\right> \geq 0$ for all $v \in K(w)$ if and only if $w=P_{K(w)}\left[w-\varrho F(w)\right]$, where $\varrho>0$ is a constant.
\end{lemma}

\begin{remark}\label{R-QVI}
It is worth noting that if $\Omega_0 \subset \mathcal{H}$ and $K:\Omega_0 \rightarrow 2^{\Omega_0}$ has convex and closed values and $F:\Omega_0 \rightarrow \mathcal{H}$ in Lemma \ref{QVI}, then the equivalent result in Lemma \ref{QVI} still holds.
\end{remark}

\begin{lemma}\label{ML-ineq}\cite{YGD} Let $k\geq 0$ and $q> 0$ be two constants, and $w$ be a locally integrable function with nonnegative valued on $[0,T)$ with $T\leq +\infty$. Assume that $z:[0,T)\to R_+$ is  locally integrable with
$$z(s) \leq w(s) +k \int_{0}^{s}(s-\varsigma)^{q-1}z(\varsigma))d\varsigma.$$
Then
$$z(s)\leq w(s) + \int_{0}^{s}  \sum\limits_{n=1}^\infty \frac{(k \Gamma(q))^n}{\Gamma(nq)} (s-\varsigma)^{nq-1} w(\varsigma) d\varsigma, \quad s\in [0,T).$$
Furthermore, if $w$ is a nondecreasing function, then $z(s)\leq w(s)  \mathrm{E}_q\left(k\Gamma(q) s^q\right)$ for all $s\in [0,T)$.
\end{lemma}

\section{The unique solvability of NFDQVI}
\noindent \setcounter{equation}{0}
In this section, the unique solvability of NFDQVI (\ref{GFDQVI}) is proved by Banach fixed point principle. In the sequel, we adopt certain assumptions regarding the data:
\begin{enumerate}
  \item[(H$_1$)] there exist two constants $l_G, \eta_G>0$ such that $G: I\times \mathcal{H}_1 \times \Omega \rightarrow \mathcal{H}_2$ is $l_G$-Lipschitz and strongly pseudomonotone w.r.t. its third argument on $K(u)$ for all $u \in \Omega$, that is,
   $$\|G(s_2,x_2,u_2)-G(s_1,x_1,u_1)\| \leq l_G (|s_2-s_1|+\|x_2-x_1\|+\|u_2-u_1\|)$$
   for all $s_1,s_2 \in I, x_1, x_2 \in \mathcal{H}_1,u_1, u_2 \in \Omega$, and
   $$\left<G(s,x,u_2),u_1-u_2\right> \geq 0 \;\;\Rightarrow \;\; \left<G(s,x,u_1),u_1-u_2\right> \geq \eta_G \|u_1-u_2\|^2$$
   for all $s \in I, x \in \mathcal{H}_1, u_1, u_2 \in K(u)$ and $u\in \Omega$;
  \item[(H$_2$)] there is a constant $l_K\geq0$ such that
   $$\big\|P_{K(u_2)}\left[\varpi\right] - P_{K(u_1)}\left[\varpi\right]\big\| \leq l_K \|u_2-u_1\|$$
   for all $u_1, u_2 \in \Omega$ and  $\varpi \in \mathcal{H}_2$;
  \item[(H$_3$)] $g: I \times \mathcal{H}_1 \times \mathcal{H}_2 \rightarrow \mathcal{H}_1$ is a continuous mappings, and there is $l_g>0$ satisfying
      $$\|g(s,x_2,u_2)-g(s,x_1,u_1)\| \leq l_g(\|x_2-x_1\|+\|u_2-u_1\|)$$
      for all $s \in I,\; x_1, x_2 \in \mathcal{H}_1$ and $u_1, u_2 \in \mathcal{H}_2$;
  \item[(H$_4$)] $f: I\times \mathcal{H}_1 \rightarrow \mathcal{H}_1$ is a continuous mappings, and there is $l_f>0$ satisfying
      $$\|f(s,x_2)-f(s,x_1)\| \leq l_f \|x_2-x_1\|,\; \forall  \, s \in I,\; x_1, x_2 \in \mathcal{H}_1;$$
   \item[(H$_5$)] $\psi: C(I,\mathcal{H}_1)  \rightarrow \mathcal{H}_1$ is a continuous mappings, and there is $0<l_\psi<1$  satisfying
      $$\|\psi(y_2)-\psi(y_1)\| \leq l_\psi \|y_2-y_1\|_\mathcal{B},\; \forall  \,  y_1, y_2 \in C(I,\mathcal{H}_1).$$
\end{enumerate}

\begin{remark}
Hypothesis (H$_2$) can be easily verified for some certain situations. For example, (i) if $\Omega = \mathcal{H}_2$, $K(u)=\phi(u)+\Omega_0$ with $\phi: \mathcal{H}_2 \rightarrow \mathcal{H}_2$ and $\Omega_0$ being a nonempty convex and closed set of $\mathcal{H}_2$, then
\begin{equation}\label{PKOmega}
P_{K(u)}\left[\varpi\right]=P_{\phi(u)+\Omega_0}\left[\varpi\right]=\phi(u) + P_{\Omega_0}\left[\varpi -\phi(u)\right]
\end{equation}
for all $u \in \mathcal{H}_2$ and $\varpi \in \mathcal{H}_2$. In addition, for any $u_1, u_2\in \mathcal{H}_2$ and $ \varpi \in \mathcal{H}_2$, one has
\begin{eqnarray*}
&& \big\|P_{K(u_2)}\left[\varpi\right] - P_{K(u_1)}\left[\varpi\right]\big\|
= \big\|P_{\phi(u_2)+\Omega_0}\left[\varpi\right] - P_{\phi(u_1)+\Omega_0}\left[\varpi\right]\big\| \\
&=& \left\|\big(\phi(u_2) + P_{\Omega_0}\left[\varpi - \phi(u_2)\right] \big) - \big(\phi(u_1) + P_{\Omega_0}\left[\varpi - \phi(u_1)\right] \big) \right\|  \\
&=& \left\|\big(\varpi -\phi(u_2) - P_{\Omega_0}\left[\varpi - \phi(u_2)\right] \big) - \big(\varpi-\phi(u_1) - P_{\Omega_0}\left[\varpi - \phi(u_1)\right] \big) \right\|.
\end{eqnarray*}
Let $\phi$ be $l_\phi$-Lipschitz. Using the fact
$$\big\|\varpi_2-P_{\Omega_0}\left[\varpi_2\right]-(\varpi_1-P_{\Omega_0}\left[\varpi_1\right])\big\| \leq \|\varpi_2-\varpi_1\|$$
for all $\varpi_1,\varpi_2 \in \mathcal{H}_2$ (see, e.g., \cite{LLH}), one has
\begin{equation}\label{u-PKOmega}
\big\|P_{K(u_2)}\left[\varpi\right] - P_{K(u_1)}\left[\varpi\right]\big\| = \big\|P_{\phi(u_2)+\Omega_0}\left[\varpi\right] - P_{\phi(u_1)+\Omega_0}\left[\varpi\right]\big\| \leq \|\phi(u_2) - \phi(u_1)\| \leq l_\phi \|u_2-u_1\|
\end{equation}
for all $u_1, u_2 \in \mathcal{H}_2$ and $ \varpi \in \mathcal{H}_2$, that is, hypothesis (H$_2$) holds. (ii) Let $K(u)=\phi(u)+\Omega$ with $\phi: \Omega \rightarrow \Omega$ and $\Omega$ being a convex and closed cone of $\mathcal{H}_2$. Similar to the argument of (\ref{u-PKOmega}), we have that hypothesis (H$_2$) also holds.
\end{remark}

Now we will show some properties of solution for PQVI in (\ref{GFDQVI})  under the hypotheses of strong pseudomonotonicity and Lipschitzean.
\begin{lemma}\label{PQVI-u0}
Let (H$_1$) and (H$_2$) hold. If
\begin{equation}\label{KlG-1}
 2 l_K^2 \left(l_G \sqrt{l_G^2+\eta_G^2} + l_G^2+\eta_G^2 \right) <\eta_G^2,
\end{equation}
then for fixed $s\in I$ and $x \in C(I,\mathcal{H}_1)$, there exists a unique solution $u(s) \in K(u(s)) \subset \Omega$ solving PQVI in (\ref{GFDQVI}).

\textbf{Proof.}\hspace{0.2cm} Fix $s\in I$ and $x \in C(I,\mathcal{H}_1)$. Let $u:I \rightarrow \Omega$ be a given function. Set $\hat{u}=u(s)$ and $\hat{x}=x(s)$. Fixing $\hat{u} \in \Omega$, similar to the argument of \cite[Theorem 3.1]{NQ}, we first consider an auxiliary problem: find $\hat{u}^* \in K(\hat{u}) \subset \Omega$ such that
\begin{equation}\label{AQVI}
  \left<G\left(s,\hat{x},\hat{u}^*\right),v-\hat{u}^*\right>\geq0, \; \forall \, v \in K(\hat{u}).
\end{equation}
It follows from hypothesis (H$_1$) and \cite[Theorem 2.1]{KVK} that the inequality (\ref{AQVI}) is unique solvable. Define a mapping $\Phi: \Omega \rightarrow \Omega$  as follows
$$\Phi(\hat{u})=\hat{u}^*, \quad  \forall \, \hat{u} \in \Omega,$$
where $\hat{u}^*$ is the unique solution of (\ref{AQVI}), namely, for any $\hat{u} \in \Omega$, $\Phi(\hat{u})$ is unique in $K\left(\hat{u}\right)$ such that
\begin{equation}\label{AQVI1} \left<G\left(s,\hat{x},\Phi(\hat{u})\right),v-\Phi(\hat{u})\right>\geq0, \; \forall \, v \in K(\hat{u}).
\end{equation}
Clearly, the unique solvable of PQVI in (\ref{GFDQVI}) is equivalent to $\Phi$ has a unique fixed point.

Next, we prove that $\Phi$ is a contractive mapping. Let $\varrho_1 >0$ be a constant. For any $\hat{u}_1, \hat{u}_2\in  \Omega$, in light of Lemma \ref{POmega} and (\ref{AQVI1}), one has
$$\Phi(\hat{u}_1)=P_{K(\hat{u}_1)}\left[\Phi(\hat{u}_1)-\varrho_1 G(s,\hat{x},\Phi(\hat{u}_1))\right]$$
and
$$\Phi(\hat{u}_2)=P_{K(\hat{u}_2)}\left[\Phi(\hat{u}_2)-\varrho_1 G(s,\hat{x},\Phi(\hat{u}_2))\right]$$
with $\Phi(\hat{u}_1) \in K(\hat{u}_1)$ and $\Phi(\hat{u}_2) \in K(\hat{u}_2)$. Let
$$z=P_{K(\hat{u}_2)}\left[\Phi(\hat{u}_1)-\varrho_1 G(s,\hat{x},\Phi(\hat{u}_1))\right].$$
In view of Lemma \ref{POmega}, we have
$$z \in K(\hat{u}_2) \;\; \mbox{and}\;\; \left<z-\Phi(\hat{u}_1)+\varrho_1 G(s,\hat{x},\Phi(\hat{u}_1)),v-z\right> \geq 0, \; \forall \, v \in K(\hat{u}_2).$$
Writing $v=\Phi(\hat{u}_2)$ yields
$$\left<z-\Phi(\hat{u}_1)+\varrho_1 G(s,\hat{x},\Phi(\hat{u}_1)),\Phi(\hat{u}_2)-z\right> \geq 0.$$
This implies that
\begin{eqnarray}\label{Phi-u12}
 \left<\Phi(\hat{u}_1)-z,\Phi(\hat{u}_2)-z\right>
&\leq& \varrho_1\left< G(s,\hat{x},\Phi(\hat{u}_1)),\Phi(\hat{u}_2)-z\right> \nonumber\\
&=& \varrho_1\left< G(s,\hat{x},\Phi(\hat{u}_1))-G(s,\hat{x},z),\Phi(\hat{u}_2)-z\right> - \varrho_1\left<G(s,\hat{x},z),z-\Phi(\hat{u}_2)\right>.
\end{eqnarray}
Recalled that $z \in K(\hat{u}_2)$, by virtue of (\ref{AQVI1}), one gets
\begin{equation*} \left<G\left(s,\hat{x},\Phi(\hat{u}_2)\right),z-\Phi(\hat{u}_2)\right>\geq0.
\end{equation*}
By the strong pseudomonotonicity of $G$ w.r.t. its third argument, one has
\begin{equation} \label{G-eta}
\left<G\left(s,\hat{x},z\right),z-\Phi(\hat{u}_2)\right> \geq \eta_G \|z-\Phi(\hat{u}_2)\|^2.
\end{equation}
Combining hypothesis (H$_1$), (\ref{Phi-u12}) and (\ref{G-eta}), we obtain
\begin{eqnarray}\label{2Phi-u}
&& 2\left<\Phi(\hat{u}_1)-z,\Phi(\hat{u}_2)-z\right> \nonumber\\
&\leq& 2\varrho_1\| G(s,\hat{x},\Phi(\hat{u}_1))-G(s,\hat{x},z)\|\|\Phi(\hat{u}_2)-z\| - 2\varrho_1 \eta_G \|z-\Phi(\hat{u}_2)\|^2 \nonumber\\
&\leq& 2\varrho_1 l_G\| \Phi(\hat{u}_1)-z\|\|\Phi(\hat{u}_2)-z\| - 2\varrho_1 \eta_G \|z-\Phi(\hat{u}_2)\|^2 \nonumber\\
&\leq& \varrho_1^2 l_G^2\| \Phi(\hat{u}_1)-z\|^2+\|\Phi(\hat{u}_2)-z\|^2 - 2\varrho_1 \eta_G \|z-\Phi(\hat{u}_2)\|^2 \nonumber\\
&=& \varrho_1^2 l_G^2\| \Phi(\hat{u}_1)-z\|^2+\|(\Phi(\hat{u}_2)-\Phi(\hat{u}_1)) + (\Phi(\hat{u}_1)-z)\|^2 - 2\varrho_1 \eta_G \|(\Phi(\hat{u}_2)-\Phi(\hat{u}_1))-(z-\Phi(\hat{u}_1))\|^2 \nonumber\\
&\leq& \varrho_1^2 l_G^2\| \Phi(\hat{u}_1)-z\|^2+\|\Phi(\hat{u}_2)-\Phi(\hat{u}_1)\|^2+\|\Phi(\hat{u}_1)-z\|^2 +2\left<\Phi(\hat{u}_2)-\Phi(\hat{u}_1),\Phi(\hat{u}_1)-z\right> \nonumber\\
&&- \varrho_1 \eta_G \|\Phi(\hat{u}_2)-\Phi(\hat{u}_1)\|^2 +2\varrho_1 \eta_G \|z-\Phi(\hat{u}_1)\|^2,
\end{eqnarray}
where the last inequality comes from the fact $\displaystyle -\|z_1-z_2\|^2 \leq -\frac{1}{2}\|z_1\|^2+\|z_2\|^2$ with $z_1,z_2 \in \mathcal{H}_2$. Note that
\begin{equation}\label{Phi-u1}
2\left<\Phi(\hat{u}_1)-z,\Phi(\hat{u}_2)-z\right> - 2\left<\Phi(\hat{u}_2)-\Phi(\hat{u}_1),\Phi(\hat{u}_1)-z\right>= 2 \|\Phi(\hat{u}_1)-z\|^2.
\end{equation}
In light of (\ref{2Phi-u}), (\ref{Phi-u1}) and hypothesis (H$_2$), one has
\begin{eqnarray}\label{Phi-u2}
&&(\varrho_1 \eta_G-1)\|\Phi(\hat{u}_2) -\Phi(\hat{u}_1)\|^2 \nonumber\\
&\leq& \left(\varrho_1^2 l_G^2 + 2\varrho_1 \eta_G -1\right) \|\Phi(\hat{u}_1)-z\|^2 \nonumber\\
&=& \left(\varrho_1^2 l_G^2 + 2\varrho_1 \eta_G -1\right) \big\|P_{K(\hat{u}_1)}\left[\Phi(\hat{u}_1)-\varrho_1 G(s,\hat{x},\Phi(\hat{u}_1))\right]-P_{K(\hat{u}_2)}\left[\Phi(\hat{u}_1)-\varrho_1 G(s,\hat{x},\Phi(\hat{u}_1))\right]\big\|^2 \nonumber\\
&\leq& \left(\varrho_1^2 l_G^2 + 2\varrho_1 \eta_G -1\right) l_K^2\|\hat{u}_1-\hat{u}_2\|^2.
\end{eqnarray}
Given $\displaystyle \varrho_1=\beta_1+\frac{1}{\eta_G}$ with $\beta_1 >0$, we have $\varrho_1 \eta_G-1=\beta_1 \eta_G>0$ and
\begin{equation}\label{Phi-u01}
\|\Phi(\hat{u}_2) -\Phi(\hat{u}_1)\|^2 \leq l_K^2 \left[\frac{l_G^2 \beta_1}{\eta_G} + \frac{1}{\beta_1}\left(\frac{l_G^2}{\eta_G^3}+\frac{1}{\eta_G}\right) + \frac{2l_G^2}{\eta_G^2}+2\right]\|\hat{u}_2-\hat{u}_1\|^2.
\end{equation}
Choosing $\beta_1$ such that
$$\frac{l_G^2 \beta_1}{\eta_G} = \frac{1}{\beta_1}\left(\frac{l_G^2}{\eta_G^3}+\frac{1}{\eta_G}\right),$$
that is, taking $\displaystyle \beta_1 = \frac{\sqrt{l_G^2+\eta_G^2}}{l_G \eta_G}$, in view of (\ref{Phi-u01}), one has
\begin{equation*}
\|\Phi(\hat{u}_2) -\Phi(\hat{u}_1)\|^2 \leq 2 \frac{l_K^2}{\eta_G^2} \left(l_G \sqrt{l_G^2+\eta_G^2} + l_G^2+\eta_G^2 \right)\|\hat{u}_2-\hat{u}_1\|^2.
\end{equation*}
It follows from (\ref{KlG-1}) that   $\displaystyle \frac{l_K}{\eta_G} \sqrt{2 \left(l_G \sqrt{l_G^2+\eta_G^2} + l_G^2+\eta_G^2 \right)} <1$, and so $\Phi$ is a contraction mapping, i.e., $\Phi$ has a unique fixed point. Consequently, for each $s \in I$ and $x \in C(I,\mathcal{H}_1)$, there exists a unique solution $u(s) \in K(u(s)) \subset \Omega$ to PQVI in (\ref{GFDQVI}). \hfill $\Box$
\end{lemma}

\begin{remark}
In this paper, we prove that PQVI in (\ref{GFDQVI}) is unique solvable under the strong pseudomonotonicity and Lipschitz condition. In 2023, Du et al. investigated the unique solvability for PQVI involved fractional differential quasi-variational inequality with Katugampola fractional operator under the strong monotonicity and Lipschitz condition (see \cite[Lemma 3.1]{DLJX}). Obviously, Lemma \ref{PQVI-u0} extends the result of \cite[Lemma 3.1]{DLJX}.
\end{remark}

\begin{lemma}\label{PQVI-u}
Let (H$_1$) and (H$_2$) hold. If it holds
\begin{equation}\label{KlG-2}
L= 2 \frac{l^2_K}{\eta^2_G} \left( l_G \sqrt{4l_G^2+2\eta_G^2} + 2l_G^2+\eta_G^2 \right) <1.
\end{equation}
Then for fixed $s\in I$ and $x \in C(I,\mathcal{H}_1)$, there exists a unique $u(s) \in K(u(s)) \subset \Omega$ solving PQVI in (\ref{GFDQVI}).
Moreover, the solution function $u:I \rightarrow \Omega$ is continuous. In addition, let $u_1 $ and $u_2 $ be two solutions to PQVI with $x$ being $x_1, x_2 \in C(I;\mathcal{H}_1)$ in (\ref{GFDQVI}), respectively. Then it holds
\begin{equation}\label{LFQ}
\|u_1(s)-u_2(s)\| \leq \sqrt{\frac{\kappa}{1-L}}\|x_1(s)-x_2(s)\|, \; \forall  \, s \in I,
\end{equation}
where $\kappa$ is defined by (\ref{kappa}).

\textbf{Proof.}\hspace{0.2cm} Let $x \in C(I,\mathcal{H}_1)$. The proof proceeds in three steps.

\textbf{Step 1.} We show that, for any given $s\in I$, PQVI in (\ref{GFDQVI}) has a unique solution.

Obviously, (\ref{KlG-2}) implies (\ref{KlG-1}). It follows from Lemma \ref{PQVI-u0} that, for any given $s\in I$, $x \in C(I,\mathcal{H}_1)$, there exists a unique solution $u(s) \in K(u(s)) \subset \Omega$ solving PQVI in (\ref{GFDQVI}).

\textbf{Step 2.} We show that  the continuity of solution  $u:I \rightarrow \Omega$.

For any given $s_1,s_2 \in I$, in view of Remark \ref{R-QVI}, one has $u(s_i) \in K(u(s_i))\;(i=1,2)$ such that
\begin{equation*}
u(s_i)=P_{K(u(s_i))}\left[u(s_i)-\varrho_2 G(s_i,x(s_i),u(s_i))\right],
\end{equation*}
where $\varrho_2 >0$ is a constant. Let
$$\widehat{z}=P_{K(u(s_2))}\left[u(s_1)-\varrho_2 G(s_1,x(s_1),u(s_1))\right].$$
Using Lemma \ref{POmega}, one gets
$$\widehat{z} \in K(u(s_2)) \;\; \mbox{and}\;\; \left<\widehat{z}-u(s_1)+\varrho_2 G(s_1,x(s_1),u(s_1)),v-\widehat{z}\right> \geq 0, \; \forall \, v \in K(u(s_2)).$$
Letting $v=u(s_2)$ yields
$$\left<\widehat{z}-u(s_1)+\varrho_2 G(s_1,x(s_1),u(s_1)),u(s_2)-\widehat{z}\right> \geq 0.$$
It follows that
\begin{eqnarray}\label{Phi-u12-2}
&& \left<u(s_1)-\widehat{z},u(s_2)-\widehat{z}\right> \nonumber\\
&\leq& \varrho_2\left< G(s_1,x(s_1),u(s_1)),u(s_2)-\widehat{z}\right> \nonumber\\
&=& \varrho_2\left< G(s_1,x(s_1),u(s_1))-G(s_2,x(s_2),\widehat{z}), u(s_2)-\widehat{z}\right> - \varrho_2\left<G(s_2,x(s_2),\widehat{z}),\widehat{z}-u(s_2)\right>.
\end{eqnarray}
Noting that $\widehat{z} \in K(u(s_2))$, we have $\left<G(s_2,x(s_2),u(s_2)),\widehat{z}-u(s_2)\right>\geq0$.
It follows from hypothesis (H$_1$) that
\begin{equation} \label{G-eta-2}
\left<G(s_2,x(s_2),\widehat{z}),\widehat{z}-u(s_2)\right> \geq \eta_G \|\widehat{z}-u(s_2)\|^2.
\end{equation}
Similar to the argument of (\ref{2Phi-u}) and (\ref{Phi-u2}), it follows from (\ref{Phi-u12-2}), (\ref{G-eta-2}) and hypotheses (H$_1$)-(H$_2$) that
\begin{eqnarray}\label{3Phi-u}
&& 2\left<u(s_1)-\widehat{z},u(s_2)-\widehat{z}\right> \nonumber\\
&\leq& 2\varrho_2\| G(s_1,x(s_1),u(s_1))-G(s_2,x(s_2),\widehat{z})\|\|u(s_2)-\widehat{z}\| - 2\varrho_2 \eta_G \|\widehat{z}-u(s_2)\|^2 \nonumber\\
&\leq& 2\varrho_2 l_G(|s_1-s_2|+\|x(s_1)-x(s_2)\|)\|u(s_2)-\widehat{z}\| + 2\varrho_2 l_G\|u(s_1)-\widehat{z}\|\|u(s_2)-\widehat{z}\| - 2\varrho_2 \eta_G \|\widehat{z}-u(s_2)\|^2 \nonumber\\
&\leq& 2\varrho_2^2 l_G^2(|s_1-s_2|+\|x(s_1)-x(s_2)\|)^2+\frac{1}{2}\|u(s_2)-\widehat{z}\|^2 +2\varrho_2^2 l_G^2\|u(s_1)-\widehat{z}\|^2+\frac{1}{2}\|u(s_2)-\widehat{z}\|^2 \nonumber\\
&& - 2\varrho_2 \eta_G \|\widehat{z}-u(s_2)\|^2 \nonumber\\
&=& 2\varrho_2^2 l_G^2(|s_1-s_2|+\|x(s_1)-x(s_2)\|)^2+2\varrho_2^2 l_G^2\|u(s_1)-\widehat{z}\|^2 + \|u(s_2)-\widehat{z}\|^2 - 2\varrho_2 \eta_G \|\widehat{z}-u(s_2)\|^2 \nonumber\\
&\leq& 2\varrho_2^2 l_G^2(|s_1-s_2|+\|x(s_1)-x(s_2)\|)^2+2\varrho_2^2 l_G^2\|u(s_1)-\widehat{z}\|^2 + \|u(s_2)-u(s_1)\|^2 + \|u(s_1)-\widehat{z}\|^2\nonumber\\
&&+2\left<u(s_2)-u(s_1),u(s_1)-\widehat{z}\right>- \varrho_2 \eta_G \|u(s_2)-u(s_1)\|^2 +2\varrho_2 \eta_G \|\widehat{z}-u(s_1)\|^2
\end{eqnarray}
and
\begin{eqnarray}\label{Phi-u3}
&&(\varrho_2 \eta_G-1)\|u(s_2)-u(s_1)\|^2 \nonumber\\
&\leq& 2\varrho_2^2 l_G^2(|s_1-s_2|+\|x(s_1)-x(s_2)\|)^2 + \left(2\varrho_2^2 l_G^2 + 2\varrho_2 \eta_G -1\right) \|u(s_1)-\widehat{z}\|^2 \nonumber\\
&=& 2\varrho_2^2 l_G^2(|s_1-s_2|+\|x(s_1)-x(s_2)\|)^2 + \left(2\varrho_2^2 l_G^2 + 2\varrho_2 \eta_G -1\right) \big\|P_{K(u(s_1))}\left[u(s_1)-\varrho_2 G(s_1,x(s_1),u(s_1))\right]\nonumber\\
&&-P_{K(u(s_2))}\left[u(s_1)-\varrho_2 G(s_1,x(s_1),u(s_1))\right]\big\|^2 \nonumber\\
&\leq& 2\varrho_2^2 l_G^2(|s_1-s_2|+\|x(s_1)-x(s_2)\|)^2 + \left(2\varrho_2^2 l_G^2 + 2\varrho_2 \eta_G -1\right) l_K^2\|u(s_2)-u(s_1)\|^2.
\end{eqnarray}
Letting $\displaystyle \varrho_2=\beta_2+\frac{1}{\eta_G}$ with $\beta_2 >0$, we have $\varrho_2 \eta_G-1=\beta_2 \eta_G>0$ and
\begin{eqnarray}\label{Phi-u01-2}
\|u(s_2)-u(s_1)\|^2
&\leq&2\frac{\left(\beta_2+\frac{1}{\eta_G}\right)^2 l_G^2}{\beta_2 \eta_G}(|s_1-s_2|+\|x(s_1)-x(s_2)\|)^2 \nonumber\\
&&+ l_K^2 \left[\frac{2 l_G^2 \beta_2}{\eta_G} + \frac{1}{\beta_2}\left(\frac{2l_G^2}{\eta_G^3}+\frac{1}{\eta_G}\right) + \frac{4 l_G^2}{\eta_G^2}+2\right]\|u(s_2)-u(s_1)\|^2.
\end{eqnarray}
Choosing $\beta_2$ such that
$$\frac{2l_G^2 \beta_2}{\eta_G} = \frac{1}{\beta_2}\left(\frac{2l_G^2}{\eta_G^3}+\frac{1}{\eta_G}\right),$$
that is, taking $\displaystyle \beta_2 = \frac{\sqrt{l_G^2+\frac{1}{2}\eta_G^2}}{l_G \eta_G}$, in view of (\ref{Phi-u01-2}), one has
\begin{equation*}
\|u(s_2)-u(s_1)\|^2 \leq  \kappa(|s_2-s_1|+\|x(s_2)-x(s_1)\|)^2+2 \frac{l_K^2}{\eta_G^2} \left(2l_G \sqrt{l_G^2+\frac{1}{2}\eta_G^2} + 2l_G^2+\eta_G^2 \right)\|u(s_2)-u(s_1)\|^2,
\end{equation*}
where
\begin{equation}\label{kappa}
\kappa=2\frac{l_G \left(\sqrt{l_G^2+\frac{1}{2}\eta_G^2}+l_G\right)^2 }{\eta_G^2\sqrt{l_G^2+\frac{1}{2}\eta_G^2}}.
\end{equation}
It follows from (\ref{KlG-2}) that
\begin{equation*}
\|u(s_2)-u(s_1)\|^2 \leq  \frac{\kappa}{1-L} (|s_2-s_1|+\|x(s_2)-x(s_1)\|)^2.
\end{equation*}
From the continuity of $x$, one gets the continuity of  $u:I \rightarrow \Omega$.

\textbf{Step 3.} We prove that the inequality (\ref{LFQ}) holds.

For any $s \in I$, let $u_1(s)$ and $u_2(s)$ be the unique solution to PQVI in (\ref{GFDQVI}) with $x_1, x_2 \in C(I;\mathcal{H}_1)$, respectively. Then by setting $\hat{u}_i=u_i(s)$ and $\hat{x}_i=x_i(s)$ with $i=1,2$, Remark \ref{R-QVI} yields that
$$
\hat{u}_i=P_{K(\hat{u}_i)}\left[\hat{u}_i-\varrho_2 G(s,\hat{x}_i,\hat{u}_i)\right],
$$
where $\varrho_2$ is defined in step 2. Let $\widetilde{z}=P_{K(\hat{u}_2)}\left[\hat{u}_1-\varrho_2 G(s,\hat{x}_1,\hat{u}_1)\right]$. Similar to the argument of (\ref{3Phi-u}) and (\ref{Phi-u3}), one has
\begin{eqnarray*}
 2\left<\hat{u}_1-\widetilde{z},\hat{u}_2-\widetilde{z}\right> &\leq& 2\varrho_2\|G(s,\hat{x}_1,\hat{u}_1)- G(s,\hat{x}_2,\widetilde{z})\|\|\hat{u}_2-\widetilde{z}\| - 2\varrho_2 \eta_G \|\widetilde{z}-\hat{u}_2\|^2 \nonumber\\
&\leq& 2\varrho_2 l_G\|\hat{x}_1-\hat{x}_2\|\|\hat{u}_2-\widetilde{z}\| + 2\varrho_2 l_G\|\hat{u}_1-\widetilde{z}\|\|\hat{u}_2-\widetilde{z}\| - 2\varrho_2 \eta_G \|\widetilde{z}-\hat{u}_2\|^2 \nonumber\\
&\leq& 2\varrho_2^2 l_G^2\|\hat{x}_1-\hat{x}_2\|^2+2\varrho_2^2 l_G^2\|\hat{u}_1-\widetilde{z}\|^2 + \|\hat{u}_2-\hat{u}_1\|^2 + \|\hat{u}_1-\widetilde{z}\|^2\nonumber\\
&&\mbox{} +2\left<\hat{u}_2-\hat{u}_1,\hat{u}_1-\widetilde{z}\right> - \varrho_2 \eta_G \|\hat{u}_2-\hat{u}_1\|^2 +2\varrho_2 \eta_G \|\widetilde{z}-\hat{u}_1\|^2
\end{eqnarray*}
and
$$(\varrho_2 \eta_G-1)\|\hat{u}_2-\hat{u}_1\|^2
\leq 2\varrho_2^2 l_G^2\|\hat{x}_1-\hat{x}_2\|^2 + \left(2\varrho_2^2 l_G^2 + 2\varrho_2 \eta_G -1\right) l_K^2\|\hat{u}_2-\hat{u}_1\|^2.
$$
It follows that
\begin{equation*}
\|\hat{u}_2-\hat{u}_1\|^2 \leq  \frac{\kappa}{1-L} \|\hat{x}_1-\hat{x}_2\|^2,
\end{equation*}
where $L$ and $\kappa$ are defined by (\ref{KlG-2}) and (\ref{kappa}), respectively. Hence the inequality (\ref{LFQ}) holds.  \hfill $\Box$
\end{lemma}

\begin{remark}\label{FDQVIx0}
If all hypotheses of Lemma \ref{PQVI-u} are satisfied, then for any given $x \in C(I,\mathcal{H}_1)$, PQVI in (\ref{GFDQVI}) has a unique solution $u_x \in C(I,\Omega)$. Define a mapping $\Lambda:C(I,\mathcal{H}_1) \rightarrow C(I,\Omega)$ by setting
\begin{equation}\label{Lambdax}
\Lambda(x)=u_x, \;\; \forall \; x \in C(I,\mathcal{H}_1)
\end{equation}
with $u_x \in C(I,\Omega)$ being the unique solution to PQVI in (\ref{GFDQVI}) associated to $x \in C(I,\mathcal{H}_1)$. We conclude from Definition \ref{def-FDQVI} that the trajectory $x$ of NFDQVI (\ref{GFDQVI}) can be given by
\begin{equation}\label{FDE}
x(s) = x_{0} + \psi(x) +\frac{1}{\Gamma(q)} \int_0^s(s-\zeta)^{q-1}\left[f(\zeta,x(\zeta)) + g(\zeta,x(\zeta),(\Lambda x)(\zeta))\right]d\zeta , \; s \in I.
\end{equation}
\end{remark}

Next we give our main result on the unique solvability of NFDQVI (\ref{GFDQVI}) as follows.

\begin{theorem}\label{exist-GFDQVI}
Let (H$_1$)-(H$_5$) and (\ref{KlG-2}) hold. Then NFDQVI (\ref{GFDQVI}) is unique solvable.

\textbf{Proof.}\hspace{0.2cm}  Consider the mapping $\Psi:C(I,\mathcal{H}_1) \rightarrow C(I,\mathcal{H}_1)$ given by
\begin{equation*}\label{Psi}
(\Psi x)(s)= x_{0} + \psi(x) +\frac{1}{\Gamma(q)} \int_0^s(s-\zeta)^{q-1}\left[f(\zeta,x(\zeta)) + g(\zeta,x(\zeta),(\Lambda x)(\zeta))\right]d\zeta , \; s \in I.
\end{equation*}
We now show that there exists a unique fixed point of $\Psi$. In fact, for any $x,y \in C(I,\mathcal{H}_1)$, we have
\begin{eqnarray}\label{Psi-xy}
&&\|(\Psi x)(s)-(\Psi y)(s)\| \nonumber\\
&\leq& \frac{1}{\Gamma(q)} \int_0^s(s-\zeta)^{q-1}\left(\|f(\zeta,x(\zeta))- f(\zeta,y(\zeta))\|+ \|g(\zeta,x(\zeta),(\Lambda x)(\zeta)) - g(\zeta,y(\zeta),(\Lambda y)(\zeta))\|\right) d\zeta \nonumber\\
&&+\|\psi(x)- \psi(y)\|
\end{eqnarray}
for all $s \in I$. By (\ref{Lambdax}), hypothesis (H$_3$) and Lemma \ref{PQVI-u}, one has
\begin{eqnarray}\label{g-Lambda}
&&\|g(\zeta,x(\zeta),(\Lambda x)(\zeta)) - g(\zeta,y(\zeta),(\Lambda y)(\zeta)) \| \nonumber\\
&\leq&  l_g (\|x(\zeta) - y(\zeta)\| + \| (\Lambda x)(\zeta) - (\Lambda y)(\zeta)\| ) \nonumber\\
&\leq&  l_g\left( 1+ \sqrt{\frac{\kappa}{1-L}}\right)\|x(\zeta) - y(\zeta)\| =\xi \|x(\zeta) - y(\zeta)\|,
\end{eqnarray}
where $\displaystyle \xi = l_g\left( 1+ \sqrt{\frac{\kappa}{1-L}}\right)$. Combining hypotheses (H$_4$)-(H$_5$), (\ref{Psi-xy}) and (\ref{g-Lambda}), we obtain
\begin{eqnarray*}
&&e^{-\gamma s}\|(\Psi x)(s)-(\Psi y)(s)\| \nonumber\\
&\leq& l_\psi e^{-\gamma s}\|x - y\|_\mathcal{B} + \frac{1}{\Gamma(q)} \int_0^s(s-\zeta)^{q-1}(l_f+ \xi)e^{-\gamma s}\|x(\zeta) - y(\zeta)\| d\zeta \nonumber\\
&\leq& l_\psi \|x - y\|_\mathcal{B} + \frac{l_f+ \xi}{\Gamma(q)} \|x - y\|_\mathcal{B} \int_0^s(s-\zeta)^{q-1}e^{-\gamma (s-\zeta)} d\zeta  \nonumber\\
&\leq& \left(l_\psi + \frac{l_f+ \xi}{\Gamma(q)} \underset{s\in I}\sup \int_0^s(s-\zeta)^{q-1}e^{-\gamma (s-\zeta)} d\zeta \right)\|x - y\|_\mathcal{B} \nonumber\\
&=& \lambda \|x - y\|_\mathcal{B},
\end{eqnarray*}
where
$$\lambda = l_\psi + \frac{l_f+ \xi}{\Gamma(q)} \underset{s\in I}\sup \int_0^s(s-\zeta)^{q-1}e^{-\gamma (s-\zeta)} d\zeta.$$
Furthermore, we take $\gamma>0$ is large enough such that
\begin{equation}\label{B-q}
 \frac{l_f+ \xi}{\Gamma(q)} \underset{s\in I}\sup \int_0^s(s-\zeta)^{q-1}e^{-\gamma (s-\zeta)} d\zeta <1 - l_\psi.
\end{equation}
This implies that $\lambda<1$. In fact, since
$$\int_0^s (s-\zeta)^{q-1} e^{-\gamma(s-\zeta)} d\zeta=\int_0^s \tau^{q-1} e^{-\gamma\tau} d \tau =\frac{1}{\gamma^q} \int_0^{\gamma s} \varsigma^{q-1} e^{-\varsigma} d\varsigma \leq \frac{1}{\gamma^q} \int_0^{+\infty} \varsigma^{q-1} e^{-\varsigma} d\varsigma = \frac{\Gamma(q)}{\gamma^q},$$
combining the condition $l_\psi <1$ in hypothesis (H$_5$), there exists sufficiently enough $\gamma>0$ such that (\ref{B-q}) holds. It follows that
$$\|\Psi x-\Psi y\|_\mathcal{B}  \leq \lambda \|x-y\|_\mathcal{B},\; \forall \, x,y \in C(I,\mathcal{H}_1) \quad \mbox{and} \quad \lambda<1,$$
that is, $\Psi$ is contractible. By using Banach fixed point principle, one gets that there exists a unique fixed point $\bar{x}$ of $\Psi$ in $C(I,\mathcal{H}_1)$ and consequently (\ref{FDE}) has a unique solution. Let $\bar{u}=\Lambda\left(\bar{x}\right)$. By Remark \ref{FDQVIx0}, it follows that $(\bar{x},\bar{u}) \in C(I,\mathcal{H}_1) \times C(I,\Omega)$ is the unique solution of NFDQVI (\ref{GFDQVI}).  \hfill $\Box$
\end{theorem}

\begin{corollary}\label{exist-GFDQCP}
Let $K(u)=\phi(u)+\Omega$ with $\phi: \Omega \rightarrow \Omega$ and $\Omega$ being a convex and closed cone of $\mathcal{H}_2$. Let (H$_1$)-(H$_5$) and (\ref{KlG-2}) hold. Then NFDQCP (\ref{GFDQCP}) is unique solvable.

\textbf{Proof.}\hspace{0.2cm} Applying Theorem \ref{exist-GFDQVI} and the equivalence of PQVI (\ref{PQVI}) and PQCP (\ref{PQCP}), it follows immediately that the unique solvability for NFDQCP (\ref{GFDQCP}).
\end{corollary}

\section{The Mittag-Leffler-Hyers-Ulam stability of NFDQVI}
\noindent \setcounter{equation}{0}
In this section, our aim is to study the MLHU stability of NFDQVI (\ref{GFDQVI}).  To this end, we give the following  result first.

\begin{lemma}\label{y-ineq}
Let (H$_1$)-(H$_5$) and (\ref{KlG-2}) hold. Assume that $(z,\upsilon) \in C(I,\mathcal{H}_1) \times C(I,\Omega)$ is a solution of inequality system (\ref{MLUH}). Then, for all $s \in I$, it holds
\begin{equation*}\label{yFDE}
\left\|z(s)-z(0)  - \frac{1}{\Gamma(q)} \int_0^s(s-\zeta)^{q-1}\left[f(\zeta,z(\zeta)) + g(\zeta,z(\zeta),(\Lambda z)(\zeta))\right]d\zeta \right\| \leq \frac{s^q}{\Gamma(q+1)}  \varepsilon.
\end{equation*}

\textbf{Proof.}\hspace{0.2cm} Let $(z,\upsilon) \in C(I,\mathcal{H}_1) \times C(I,\Omega)$ be a solution of (\ref{MLUH}). It follows from Remark \ref{ISys} that there is a continuous function $h:I \rightarrow \mathcal{H}_1$ satisfying
\begin{itemize}
  \item[(i)] for any $s\in I$, $\|h(s)\|\leq \varepsilon$;
  \item[(ii)] $(z,\upsilon)$ is the solution of FDQVI
  $$
\left\{
\begin{array}{l}
\tensor*[^{C}_0]{\mathrm{D}}{_s^q} z(s) = f(s,z(s)) + g(s,z(s),\upsilon(s)) + h(s), \; \forall  \, s \in I,\\
\upsilon(s) \in \mbox{SOL}(K\left(\cdot),G(s,z(s),\cdot)\right), \; \forall  \, s \in I, \\
z\left(0\right)=x\left(0\right).
\end{array}
\right.
$$
\end{itemize}
By Definition \ref{def-FDQVI} and Remark \ref{FDQVIx0}, we have
$$
z(s)= z\left(0\right) + \frac{1}{\Gamma(q)} \int_0^s(s-\zeta)^{q-1}\left[f(\zeta,z(\zeta)) + g(\zeta,z(\zeta),(\Lambda z)(\zeta))\right]d\zeta + \frac{1}{\Gamma(q)} \int_0^s(s-\zeta)^{q-1}h(\zeta)d\zeta
$$
and so
\begin{eqnarray}\label{y-x0}
 &&\left\|z(s)-z(0) - \frac{1}{\Gamma(q)} \int_0^s(s-\zeta)^{q-1}\left[f(\zeta,z(\zeta)) + g(\zeta,z(\zeta),(\Lambda z)(\zeta))\right]d\zeta \right\| \nonumber\\
 &=& \left\|\frac{1}{\Gamma(q)} \int_0^s(s-\zeta)^{q-1}h(\zeta)d\zeta \right\|  \nonumber\\
 &\leq&\frac{\varepsilon}{\Gamma(q)} \int_0^s(s-\zeta)^{q-1}d\zeta\nonumber\\
 &=& \frac{s^q}{\Gamma(q+1)}  \varepsilon, \; \forall  \, s \in I.
\end{eqnarray}
It yields the result.   \hfill $\Box$
\end{lemma}

Now we will establish the MLHU stability of NFDQVI (\ref{GFDQVI}).

\begin{theorem}\label{GFDQVI-MLUH}
If all hypotheses of Theorem \ref{exist-GFDQVI} hold, then NFDQVI (\ref{GFDQVI}) is MLHU stable.

\textbf{Proof.}\hspace{0.2cm} For any $\varepsilon >0$, let $(z,\upsilon) \in C(I,\mathcal{H}_1) \times C(I,\Omega)$ be a solution of (\ref{MLUH}) and $(x,u) \in C(I,\mathcal{H}_1) \times C(I,\Omega)$ be the unique solution of NFDQVI (\ref{GFDQVI}). Combining the nonlocal condition $z(0)=x(0)=x_0+\psi(x)$ and Remark \ref{FDQVIx0}, we have
\begin{eqnarray}\label{xyfg}
 &&\left\|z(s)-x(s) \right\|  \nonumber\\
 &=& \left\|z(s)-x_0-\psi(x)-\frac{1}{\Gamma(q)} \int_0^s(s-\zeta)^{q-1}\left[f(\zeta,x(\zeta)) + g(\zeta,x(\zeta),(\Lambda x)(\zeta))\right]d\zeta \right\| \nonumber\\
 &\leq&\frac{1}{\Gamma(q)} \int_0^s(s-\zeta)^{q-1}\left(\|f(\zeta,z(\zeta))-f(\zeta,x(\zeta))\| + \|g(\zeta,z(\zeta),(\Lambda z)(\zeta)) - g(\zeta,x(\zeta),(\Lambda x)(\zeta)) \|\right)d\zeta\nonumber\\
&&+\left\|z(s)-z(0) - \frac{1}{\Gamma(q)} \int_0^s(s-\zeta)^{q-1}\left[f(\zeta,z(\zeta)) + g(\zeta,z(\zeta),(\Lambda z)(\zeta))\right]d\zeta \right\|.
\end{eqnarray}
In view of Lemma \ref{y-ineq}, (\ref{g-Lambda}), (\ref{xyfg}), hypotheses (H$_3$) and (H$_4$), we obtain
\begin{eqnarray}\label{xyfg-1}
 &&\left\|z(s)-x(s) \right\|  \nonumber\\
  &\leq& \frac{s^q}{\Gamma(q+1)}  \varepsilon + \frac{l_f+ \xi}{\Gamma(q)} \int_0^s(s-\zeta)^{q-1}\|z(\zeta) - x(\zeta)\| d\zeta\nonumber\\
  &\leq& \frac{T^q}{\Gamma(q+1)}  \varepsilon + \frac{l_f+ \xi}{\Gamma(q)} \int_0^s(s-\zeta)^{q-1}\|z(\zeta) - x(\zeta)\| d\zeta,
\end{eqnarray}
where $\xi$ is defined by (\ref{g-Lambda}). Let
\begin{equation}\label{xyfg-2}
a = \frac{T^q}{\Gamma(q+1)},\quad b = \frac{l_f+ \xi}{\Gamma(q)}.
\end{equation}
Using Lemma \ref{ML-ineq}, one has
$$\left\|z(s)-x(s) \right\| \leq a \varepsilon \mathrm{E}_q(b \Gamma(q) s^q),$$
which implies that NFDQVI (\ref{GFDQVI}) is MLHU stable.  \hfill $\Box$
\end{theorem}

Similar to the argument of Corollary \ref{exist-GFDQCP} and Theorem \ref{GFDQVI-MLUH}, one has the following result for NFDQCP (\ref{GFDQCP}).

\begin{corollary}
If all hypotheses of Corollary \ref{exist-GFDQCP} hold, then NFDQCP (\ref{GFDQCP}) is MLHU stable.
\end{corollary}

Next we present the result of generalized MLHUR stability of NFDQVI (\ref{GFDQVI}).

\begin{theorem}\label{GFDQVI-GMLUHR}
Let all hypotheses of Theorem \ref{exist-GFDQVI} hold. If $\varphi \in C(I,R_+)$ is a nondecreasing function, then NFDQVI (\ref{GFDQVI}) is generalized MLHUR stable w.r.t. $\varphi \mathrm{E}_q$.

\textbf{Proof.}\hspace{0.2cm} Let $(z,\upsilon) \in C(I,\mathcal{H}_1) \times C(I,\Omega)$ solve the inequality system (\ref{GMLUHR}) and $(x,u) \in C(I,\mathcal{H}_1) \times C(I,\Omega)$ be the unique solution of NFDQVI (\ref{GFDQVI}). Since $\varphi \in C(I,R_+)$ is a nondecreasing function, similar to the proof of (\ref{y-x0}), it holds
\begin{align*}
 &\left\|z(s)-z(0)  - \frac{1}{\Gamma(q)} \int_0^s(s-\zeta)^{q-1}\left[f(\zeta,z(\zeta)) + g(\zeta,z(\zeta),(\Lambda z)(\zeta))\right]d\zeta \right\| \nonumber\\
 \leq& \frac{1}{\Gamma(q)} \int_0^s(s-\zeta)^{q-1}\varphi(\zeta)d\zeta \leq \frac{s^q}{\Gamma(q+1)}\varphi(s) \leq \frac{T^q}{\Gamma(q+1)}\varphi(s).
\end{align*}
Similar to the proof of (\ref{xyfg-1}), one has
\begin{eqnarray*}
 &&\left\|z(s)-x(s) \right\|  \nonumber\\
  &\leq& \left\|z(s)-z(0)  - \frac{1}{\Gamma(q)} \int_0^s(s-\zeta)^{q-1}\left[f(\zeta,z(\zeta)) + g(\zeta,z(\zeta),(\Lambda z)(\zeta))\right]d\zeta \right\|\nonumber\\
&&+ \frac{1}{\Gamma(q)} \int_0^s(s-\zeta)^{q-1}\left(\|f(\zeta,z(\zeta))-f(\zeta,x(\zeta))\| + \|g(\zeta,z(\zeta),(\Lambda z)(\zeta)) - g(\zeta,x(\zeta),(\Lambda x)(\zeta)) \|\right)d\zeta\nonumber\\
&\leq& \frac{T^q}{\Gamma(q+1)}\varphi(s) + \frac{l_f+ \xi}{\Gamma(q)} \int_0^s(s-\zeta)^{q-1}\|z(\zeta) - x(\zeta)\| d\zeta\\
&=& a \varphi(s) + b \int_0^s(s-\zeta)^{q-1}\|z(\zeta) - x(\zeta)\| d\zeta,
\end{eqnarray*}
where $\xi$ is the same as in (\ref{g-Lambda}), $a$ and $b$ are the same as in (\ref{xyfg-2}).
Applying Lemma \ref{ML-ineq}, we obtain
$$\left\|z(s)-x(s) \right\| \leq a \varphi(s) \mathrm{E}_q(b \Gamma(q) s^q),$$
which implies that NFDQVI (\ref{GFDQVI}) is generalized MLHUR stable  w.r.t. $\varphi \mathrm{E}_q$.  \hfill $\Box$
\end{theorem}

Similar to the argument of Corollary \ref{exist-GFDQCP} and Theorem \ref{GFDQVI-GMLUHR}, one has the following result for FDQCP (\ref{GFDQCP}).

\begin{corollary}
If all hypotheses of Corollary \ref{exist-GFDQCP} hold and $\varphi \in C(I,R_+)$ is a nondecreasing function, then NFDQCP (\ref{GFDQCP}) is generalized MLUHR stable w.r.t. $\varphi \mathrm{E}_q$.
\end{corollary}

\section{Applications}
\noindent \setcounter{equation}{0}
\vspace{-1.2cm}
\subsection{A time-dependent multi-agent optimization problem}
In this subsection, we show the unique solvability and two stability results for a time-dependent multi-agent optimization problem (MAOP for short) by employing Theorems \ref{exist-GFDQVI}, \ref{GFDQVI-MLUH} and \ref{GFDQVI-GMLUHR}.

Assume that there are $n$ agents, $x_j(s) \, (1\leq j \leq n)$ denote the state processes of agent $j$ at time $s$ ($s\in I=[0,T]$) and $u_i(s)$ denote the strategy of agent $j$ at time $s$. We will denote by $K_j(u_{-j}(s))$ the set of $j$-th agent's strategy that is dependent on the strategies of other $n-1$ agents  $$u_{-j}(s)=(u_1(s),\ldots,u_{j-1}(s),u_{j+1}(s),\ldots,u_n(s))^\top,$$
where $K_j:R^{n-1} \rightarrow 2^R$ is a set-valued mapping with nonempty convex and closed values. Let $$u(s)=\left(u_1(s),u_2(s),\ldots,u_n(s))^\top=(u_j(s),u_{-j}(s)\right)$$
and
$$x(s)=\left(x_1(s),x_2(s),\ldots,x_n(s))^\top=(x_j(s),x_{-j}(s)\right).$$
Agent $j$'s cost function $\Theta_j(x(s),u(s))$ depends on all agents' states and strategies. Following the ideas of Zeng et al. \cite{ZZH25}, we consider the multi-agent problem in the framework of a fractional differential game problem. For any $s\in I$, the general Nash equilibrium for the MAOP is to find the optimal strategy $u^*(s)$ and its corresponding optimal state $x^*(s)$ such that $$\Theta_j(x_j^*(s),x_{-j}^*(s),u_j^*(s),u_{-j}^*(s))=\underset{u_j(s) \in K_j(u_{-j}^*(s))}{\text{min}}\quad \Theta_j(x_j(s),x_{-j}^*(s),u_j(s),u^*_{-j}(s)),\; 1\leq j \leq n,$$
and the corresponding state $x_j(s)$ satisfy the following fractional order system:
\begin{equation}\label{state}
\left\{
\begin{aligned}
&\tensor*[^{C}_0]{\mathrm{D}}{_s^q} x_j(s) = f_j(s,x_j(s)) + g_j(s,x_j(s),u_j(s)),\\
&x_j(0)=x_{0j}+a_j \frac{\int_0^T x_j(\zeta)d\zeta}{T},
\end{aligned}
\right.
\end{equation}
where the initial value depends on the mean of the state process, $x_{0j}$ and $a_j$ are given real numbers. Moreover, the pair $(x^*,u^*)$ is called the optimal pair. The MAOP that we consider here can be summarized as follows:
\begin{equation}\label{multi}
\begin{aligned}
\underset{u_j(s)\in K_j(u_{-j}^*(s))}{\text{min}} & \quad \Theta_j(x_j(s),x_{-j}^*(s),u_j(s),u^*_{-j}(s)) \\
\text{subject to} & \quad \tensor*[^{C}_0]{\mathrm{D}}{_s^q} x_j(s) = f_j(s,x_j(s)) + g_j(s,x_j(s),u_j(s)), \\
& \quad x_j(0)=x_{0j}+a_j \frac{\int_0^T x_j(\zeta)d\zeta}{T}, \\
& \quad \forall \, s \in I.
\end{aligned}
\end{equation}
It is worth noting that a general Nash equilibrium for MAOP \eqref{multi} can be obtained by solving a class of NFDQVI. To this end, we assume $\Theta_j(\cdot,\cdot,\cdot,\cdot)$ is convex and continuously differentiable to its third argument, $\bar{G}(y,v)=\left(\triangledown_{v_1}\Theta_1(y,v),\triangledown_{v_2}\Theta_2(y,v), \cdots, \triangledown_{v_n}\Theta_n(y,v)\right)$ $(y=(y_1, y_2,\ldots,y_n)^\top, v=(v_1, v_2,\ldots,v_n)^\top \in R^n)$ and $K(v)=\prod_{j=1}^n K_j(v_{-j})$. Moreover, we define a function $G(s,x(s),u(s))=\bar{G}(x(s),u(s)),\forall s\in I$. Then we immediately have the following result.
\begin{theorem}\label{GNash}
The pair $(x^*,u^*)$ is a general Nash equilibrium for MAOP \eqref{multi} if and only if $(x^*,u^*)$ satisfies the following NFDQVI:
\begin{equation}\label{QVInash}
\left\{
\begin{array}{l}
\tensor*[^{C}_0]{\mathrm{D}}{_s^q} x^*(s) = f(s,x^*(s)) + g(s,x^*(s),u^*(s)), \;\forall \, s\in I, \\
 \left<G(s,x^*(s),u^*(s),v-u^*(s))\right>\geq 0,\;\forall v\in K(u^*(s)),\;\forall \, s\in I,\\
 x^*(0)=x^*_{0}+\psi(x^*),
\end{array}
\right.
\end{equation}
where $$f(s,x^*(s))=\left(f_{1}(s,x_1^*(s)),f_{2}(s,x_2^*(s)),\ldots, f_{n}(s,x^*_n(s))\right)^\top,$$ $$g(s,x^*(s),u^*(s))=\left(g_{1}(s,x_1^*(s),u_1^*(s)),g_{2}(s,x_2^*(s),u_2^*(s)),\ldots, g_{n}(s,x_n^*(s),u_n^*(s))\right)^\top,$$
\begin{equation}\label{psi12}
x^*_{0}=\left(x_{01},x_{02},\ldots, x_{0n}\right)^\top, \; \psi(x^*)=\left(\psi_1(x_1^*),\psi_2(x_2^*),\ldots, \psi_n(x_n^*)\right)^\top, \; \psi_j(x_j^*)=a_j \frac{\int_0^T x_j^*(\zeta)d\zeta}{T}.
\end{equation}

\textbf{Proof.}\hspace{0.2cm}  By convexity and minimum principle (see, e.g., \cite[Propositions 1.2 and 1.3]{Nagurney}), it follows that $\left(x^*,u^*\right)$ is the general Nash equilibrium if and only if for each $j$, $x_j^*$ satisfies the corresponding FDE, and $u_j^*$ admits
 \begin{equation}\label{QVInash1}
 \left<\triangledown_{u_j^*}\Theta_j(x_j^*(s),x_{-j}^*(s),u_j^*(s),u_{-j}^*(s)), v_j-u_j^*(s)\right>\geq 0,\;\forall v_j\in K_j(u_{-j}^*(s)),\;\forall \, s\in I.
\end{equation}
And then \eqref{QVInash} can be obtained.

Conversely, if $\left(x^*,u^*\right)$ satisfies \eqref{QVInash}, then $x^*$ satisfies the corresponding FDE and the QVI in \eqref{QVInash} is satisfied. For any fixed $j$ and $s$, in the QVI defined by \eqref{QVInash}, let $v_{-j}=u_{-j}^*(s)$ and $v_j$ is an arbitrary element in $K_j\left(u_{-j}^*(s)\right)$. Then one has \eqref{QVInash1} immediately.   \hfill $\Box$
\end{theorem}

Now we present the unique solvability and two stability results for the MAOP \eqref{multi}.

\begin{theorem}\label{MOP}
Assume that the functions in \eqref{QVInash} satisfy hypotheses (H$_1$)-(H$_4$) and (\ref{KlG-2}). If $0<\underset{j=1,\cdots,n}\max{|a_j|} < 1$, then the MAOP (\ref{multi}) is unique solvable.

\textbf{Proof.}\hspace{0.2cm} Let $\|x(s)\|=\sum\limits_{j=1}\limits^n|x_j(s)|$. For any $x^1=\left(x_1^1,x_2^1, \ldots, x_n^1\right)^\top,x^2=\left(x_1^2,x_2^2, \ldots, x_n^2\right) \in C(I,R^n)$, it follows from (\ref{psi12}) that
\begin{eqnarray*}
\left\|\psi(x^2)-\psi(x^1) \right\|
  &=& \sum\limits_{j=1}\limits^n\left|\psi_j\left(x_j^2\right)-\psi_j\left(x_j^1\right)\right|
= \sum\limits_{j=1}\limits^n\left|a_j\frac{\int_0^T x_j^2(\zeta)d\zeta}{T}-a_j\frac{\int_0^T x_j^1(\zeta)d\zeta}{T}\right|\nonumber\\
&\leq& \sum\limits_{j=1}\limits^n  \frac{|a_j|}{T}\int_0^T \left|x_j^2(\zeta) - x_j^1(\zeta)\right|d\zeta
= \sum\limits_{j=1}\limits^n  \frac{|a_j|}{T}\int_0^T e^{\gamma \zeta} e^{-\gamma \zeta} \left|x_j^2(\zeta) - x_j^1(\zeta)\right|d\zeta\\
&\leq& \sum\limits_{j=1}\limits^n  \frac{|a_j|}{T} \left|x_j^2 - x_j^1\right|_\mathcal{B} \int_0^T e^{\gamma \zeta} d\zeta
= \sum\limits_{j=1}\limits^n \frac{|a_j|}{T} \frac{e^{\gamma T}-1}{\gamma}\left|x_j^2 - x_j^1\right|_\mathcal{B}\\
&\leq& \underset{j=1,\cdots,n}\max{|a_j|}\cdot \frac{e^{\gamma T}-1}{\gamma T}\left\|x^2 - x^1\right\|_\mathcal{B}= l_\psi \left\|x^2 - x^1\right\|_\mathcal{B},
\end{eqnarray*}
where $l_\psi= \underset{j=1,\cdots,n}\max{|a_j|}\cdot \frac{e^{\gamma T}-1}{\gamma T}$. We claim that there exists $\gamma>0$ such that $0<l_\psi<1$. Indeed, it is clear that for $\gamma>0$, we have $l_\psi>0$. The remaining step is to show $l_\psi < 1$ for some $\gamma > 0$. Consider the function $\varrho: [0, +\infty) \rightarrow R$ given by
$$\varrho(\gamma)=e^{\gamma T}-\frac{T}{\underset{j=1,\cdots,n}\max{|a_j|}}\gamma-1, \; \gamma\geq0.$$
In view of $0<\underset{j=1,\cdots,n}\max{|a_j|} < 1$, it follows that $\varrho^\prime(\gamma*)=0$, $\varrho$ is decreasing on $[0, \gamma*]$ and increasing on $[\gamma*, +\infty)$, where $$\gamma*=\frac{1}{T}\ln\frac{1}{\underset{j=1,\cdots,n}\max{|a_j|}}>0.$$ Combining this with $\varrho(0)=0$, we conclude that there exists $\gamma \in (0, \gamma*]$ such that $\varrho(\gamma) <0$. Hence,
$$\underset{j=1,\cdots,n}\max{|a_j|}\cdot \frac{e^{\gamma T}-1}{\gamma T}<1,$$ that is, $l_\psi<1$. Therefore, hypothesis (H$_5$) holds. Applying Theorem \ref{exist-GFDQVI}, one has NFDQVI \eqref{QVInash} has a unique solution $(x^*,u^*)\in C(I,R^n)\times C(I,R^n)$. It follows from Theorem \ref{GNash} immediately that the MAOP (\ref{multi}) is unique solvable. \hfill $\Box$
\end{theorem}

\begin{theorem}
If all hypotheses of Theorem \ref{MOP} hold, then there exists $c>0$ such that for any $\varepsilon>0$, the optimal state $x^*$ of MAOP (\ref{multi}) satisfies
$$\|z^*(s)-x^*(s)\| \leq  c \varepsilon \mathrm{E}_q(s), \; \forall \,  s \in I,$$
where $\left(z^*,\upsilon^*\right)$ is a solution of the following inequality system
\begin{equation*}
\left\{
\begin{array}{l}
\left\|\tensor*[^{C}_0]{\mathrm{D}}{_s^q} z^*(s) -f(s,z^*(s)) - g(s,z^*(s),\upsilon^*(s))\right\| \leq \varepsilon, \;  \forall \, s \in I, \\
 \left<G(s,z^*(s),\upsilon^*(s),w-\upsilon^*(s))\right>\geq 0,\; \forall  \, w\in K(\upsilon^*(s)),\;\forall \, s\in I,\\
 z^*(0)=x^*(0).
\end{array}
\right.
\end{equation*}

\textbf{Proof.}\hspace{0.2cm} By applying Theorems \ref{GFDQVI-MLUH} and \ref{MOP}, we conclude that NFDQVI \eqref{QVInash} is MLHU stable. The result follows directly from Theorem \ref{GNash}. \hfill $\Box$
\end{theorem}

\begin{theorem}
Let all hypotheses of Theorem \ref{MOP} hold. If $\varphi \in C(I,R_+)$ is a nondecreasing function, then there exists a real number $c_\varphi>0$ such that the optimal state $x^*$ of MAOP (\ref{multi}) satisfies
$$\|z^*(s)-x^*(s)\| \leq  c_\varphi \varphi(s) \mathrm{E}_q(s), \; \forall \,  s \in I.$$
where $\left(z^*,\upsilon^*\right)$ is a solution of the following inequality system
\begin{equation*}
\left\{
\begin{array}{l}
\left\|\tensor*[^{C}_0]{\mathrm{D}}{_s^q} z^*(s) -f(s,z^*(s)) - g(s,z^*(s),\upsilon^*(s))\right\| \leq \varphi(s), \;\forall \, s \in I \\
 \left<G(s,z^*(s),\upsilon^*(s),w-\upsilon^*(s))\right>\geq 0,\; \forall  w\in K(\upsilon^*(s)), \;\forall \, s\in I,\\
 z^*(0)=x^*(0).
\end{array}
\right.
\end{equation*}

\textbf{Proof.}\hspace{0.2cm}  Combing Theorems \ref{GFDQVI-GMLUHR} and \ref{MOP}, we conclude that NFDQVI \eqref{QVInash} is generalized MLHUR stable w.r.t. $\varphi \mathrm{E}_q$. Theorem \ref{GNash} then yields the desired result. \hfill $\Box$
\end{theorem}

\subsection{A time-dependent price control problem}
In this subsection, we show the unique solvability and MLHU results for a time-dependent price control problem (PCP for short) by employing Theorems \ref{exist-GFDQVI}, \ref{GFDQVI-MLUH} and \ref{GFDQVI-GMLUHR}. In the sequel, we adopt the usual $\|\cdot\|_i$ to denote the $i$-norm in the $R^m$ space. In particular, we omit the subscripts of the  $2$-norm.

Given a market model with $m$ commodities. Let $u_j(s) (s \in I=[0,T], j= 1,2,\cdots, m)$ be the price at time $s$ for each commodity $j$. At the beginning of products sale, the manufacturer will provide a price range for each commodity $j$, i.e.,
\begin{equation}\label{ab-uj0}
a^0_j \leq u_j(0) \leq b^0_j,
\end{equation}
where $a^0_j$ and $b^0_j$ are two constants. However, the price of $j$-th commodity will also fluctuate in accordance with changes in the prices of other products. It can be assumed that
\begin{equation}\label{ab-uj}
a_j(u_{-j}(s)) \leq u_j(s) \leq b_j(u_{-j}(s)),\;\forall \, t \in(0,T],
\end{equation}
where $u_{-j}(s)=(u_1(s),\ldots,u_{j-1}(s),u_{j+1}(s),\ldots,u_m(s))^\top$, $a_j(u_{-j}(0))=a^0_j$, $b_j(u_{-j}(0))=b^0_j$. Let $u(s)=(u_1(s),u_2(s),\ldots,u_m(s))^\top=(u_j(s),u_{-j}(s))$. Following the analogous arguments in \cite{LV,JSZ}, let $S_j(u(s),x(s))$ and $D_j(u(s),x(s))$ be the supply and demand functions for $j$-th commodity at time $s\in I$, respectively, in which the outside force $x(s) \; (x(s)\in R^n)$ can be modeled by the following dynamic system:
\begin{equation}\label{GFDE0}
x^\prime(s) = \chi(s,x(s)) + \vartheta(s,x(s),u(s)), \;\forall \, s \in I,
\end{equation}
where $\chi:I \times R^n \rightarrow R^n$ and $\vartheta:I \times R^n \times R^m \rightarrow R^n$ are two given functions. Noticing that the outside force $x(s)$ can be considered as the production cost, exchange rate and so on, it is quite suitable to describe the evolutionary process of $x(s)$ by employing the fractional derivative in the economic process with memory \cite{TT,TT1}. In view of the elegant property of Caputo type fractional derivative, it would be necessary and interesting to update (\ref{GFDE0}) to the following new framework by using the fractional derivative
\begin{equation}\label{GFDE}
\tensor*[^{C}_0]{\mathrm{D}}{_s^q} x(s) = \chi(s,x(s)) + \vartheta(s,x(s),u(s)), \;\forall \, s \in I,
\end{equation}
In addition, let $0<t_1\leq t_2 \leq \cdots \leq t_k <T$. We assume that the initial value of outside force is determined by factors such as production cost and exchange rate at time $t_i \, (i=1,2, \ldots,k)$, that is,
\begin{equation}\label{iv-GFDE}
  x(0)=x_0+\sum\limits_{i=1}\limits^k \iota_i x(t_i),
\end{equation}
where $x_0$ and $\iota_i\,  (i=1,2, \ldots,k)$ are given constants. The aim of our work is to find $x : I \rightarrow R^n$ such that for any $s \in I$, the pair $(x(s), u(s))$ satisfies the market equilibrium condition, that is,
\begin{equation}\label{SD-ux}
S_j(u(s),x(s))-D_j(u(s),x(s)) \left\{\begin{array}{l}
\leq 0, \quad \mbox{if} \quad u_j(s)=b_j(u_{-j}(s)); \\
= 0, \quad \mbox{if} \quad a_j(u_{-j}(s)) < u_j(s) < b_j(u_{-j}(s)); \\
\geq 0, \quad \mbox{if} \quad u_j(s)=a_j(u_{-j}(s))
\end{array}
\right.
\end{equation}
with $j = 1,2,\ldots, m$, or equivalently
$$\left<S(u(s),x(s))-D(u(s),x(s)), v- u(s)\right> \geq 0, \; \forall \; v \in K(u(s)),$$
where
$$S(u(s),x(s))=(S_1(u(s),x(s)),S_2(u(s),x(s)),\ldots,S_m(u(s),x(s)))^\top,$$
$$D(u(s),x(s))=(D_1(u(s),x(s)),D_2(u(s),x(s)),\ldots,D_m(u(s),x(s)))^\top,$$
and
$$K(u(s))=\left\{u(s) \in R^m: a_j(u_{-j}(s)) \leq u_j(s) \leq b_j(u_{-j}(s)), j = 1,2,\ldots, m\right\}.$$
Clearly, the PCP described by (\ref{ab-uj0}), (\ref{ab-uj}), (\ref{GFDE}), (\ref{iv-GFDE}) and (\ref{SD-ux}) can be rewritten to the following coupled system:
\begin{equation}\label{GFDQVI-ex}
\left\{
\begin{array}{l}
\tensor*[^{C}_0]{\mathrm{D}}{_s^q} x(s) =  \chi(s,x(s)) + \vartheta(s,x(s),u(s)), \;\forall \, s \in I,\\
\left<S(u(s),x(s))-D(u(s),x(s)), v- u(s)\right> \geq 0, \quad \forall \; v \in K(u(s)), \;\forall \, s \in I, \\
x\left(0\right)=x_{0}+\sum\limits_{i=1}\limits^k \iota_i x(t_i).
\end{array}
\right.
\end{equation}
Obviously,  PCP (\ref{GFDQVI-ex}) is a form of NFDQVI (\ref{GFDQVI}), where $\psi(x)=\sum\limits_{i=1}\limits^k \iota_i x(t_i)$. Here and bellows, we assume that:
\begin{enumerate}
  \item[(A$_1$)] there exist constants $\eta_{SD}, l_{Sj}, l_{Dj}>0$ $(j= 1,2,\ldots, m)$ such that
   $$\sum_{j=1}^m S_j\left(u^1, x\right)\left(u_j^2-u_j^1\right) +  \sum_{j=1}^m D_j\left(u^1, x\right)\left(u_j^1-u_j^2\right)\geq0$$
   implies
   $$\sum_{j=1}^m S_j\left(u^2, x\right)\left(u_j^2-u_j^1\right) +  \sum_{j=1}^m D_j\left(u^2, x\right)\left(u_j^1-u_j^2\right) \geq \eta_{SD} \left\|u^2-u^1\right\|^2$$
   for all $x \in R^n, u^1=(u^1_1,u^1_2,\ldots, u^1_m)^\top, u^2=(u^2_1,u^2_2,\ldots, u^2_m)^\top \in R^m$, and
   $$\left|S_j(u^1, x^1)-S_j(u^2, x^2)\right| \leq l_{Sj} \left(\left\|x^1-x^2\right\|+\left\|u^1-u^2\right\|\right),$$
   $$\left|D_j(u^1, x^1)-D_j(u^2, x^2)\right| \leq l_{Dj}   \left(\left\|x^1-x^2\right\|+\left\|u^1-u^2\right\|\right)$$
   for all $x^1, x^2 \in R^n, u^1, u^2 \in R^m$, respectively;
  \item[(A$_2$)] $a_j(u_{-j})=\phi_j(u_{-j})+c_j, b_j(u_{-j})=\phi_j(u_{-j})+d_j \; (j= 1,2,\ldots, m)$, where $d_j$ and $c_j$ are positive constants with $c_j < d_j$ $(j= 1,2,\ldots, m)$, $u_{-j}=(u_1,\ldots,u_{j-1},u_{j+1},\ldots,u_m)^\top \in R^{m-1}$, $\phi_j: R^{m-1} \rightarrow R$ is $l_j$-Lipschitz such that
  $$L_1 =  2\frac{\|l\|^2_1}{\eta_{SD}^2}\left( \left(\|l_S\|_1+\|l_D\|_1\right) \sqrt{4\left(\|l_S\|_1+\|l_D\|_1\right)^2 +2\eta_{SD}^2} +2\left(\|l_S\|_1+\|l_D\|_1\right)^2 +\eta_{SD}^2\right) <1$$
  with $l=(l_1,l_2,\ldots, l_m)^\top$, $l_S=(l_{S1},l_{S2},\ldots, l_{Sm})^\top$ and $l_D=(l_{D1},l_{D2},\ldots, l_{Dm})^\top$;
  \item[(A$_3$)] $\vartheta: I \times R^n \times R^m \rightarrow R^n$ is a continuous function, and there exists a constant $l_\vartheta>0$ satisfying
      $$\left\|\vartheta\left(s,x^2,u^2\right)-\vartheta\left(s,x^1,u^1\right)\right\| \leq l_\vartheta\left(\left\|x^2-x^1\right\|+\left\|u^2-u^1\right\|\right)$$
      for all $s \in I,\; x^1, x^2 \in R^n$ and $u^1, u^2 \in R^m$;
  \item[(A$_4$)] $\chi: I\times R^n \rightarrow R^n$ is a continuous function, and there exists $l_\chi>0$ satisfying
      $$\left\|\chi\left(s,x^2\right)-\chi\left(s,x^1\right)\right\| \leq l_\chi \left\|x^2-x^1\right\|,\; \forall \, s \in I,\; x^1, x^2 \in R^n.$$
\end{enumerate}

Now we present the unique solvability and MLHU stability results for the time-dependent price control problem (\ref{GFDQVI-ex}).

\begin{theorem}\label{PCP}
Let (A$_1$)-(A$_4$) hold. If $0<\sum\limits_{i=1}\limits^k |\iota_i| <1$, then PCP (\ref{GFDQVI-ex}) is unique solvable.

\textbf{Proof.}\hspace{0.2cm} Let $G(x,u)=S(u, x)-D(u, x)$. It follows from hypothesis (A$_1$) that, for all $x^1, x^2 \in R^n,u^1, u^2 \in R^m$,
\begin{align*}
&\left\|G\left(x^2,u^2\right)-G\left(x^1,u^1\right)\right \|
\leq\left\|G\left(x^2,u^2\right)-G\left(x^1,u^1\right)\right\|_1 \\
=& \left\|S\left(u^2, x^2\right)-D\left(u^2, x^2\right) - \left(S\left(u^1, x^1\right)-D\left(u^1, x^1\right)\right)\right\|_1  \nonumber\\
\leq& \left\|S\left(u^2, x^2\right)-S\left(u^1, x^1\right) \right\|_1 +\left\|D\left(u^2, x^2\right) -D\left(u^1, x^1\right)\right\|_1 \\
=& \sum_{j=1}^m \left|S_j\left(u^2, x^2\right)-S_j\left(u^1, x^1\right)\right| +\sum_{j=1}^m  \left|D_j\left(u^2, x^2\right) -D_j\left(u^1, x^1\right)\right| \nonumber\\
\leq& \sum_{j=1}^m l_{Sj} \left(\left\|x^2-x^1\right\|+\left\|u^2-u^1\right\|\right) + \sum_{j=1}^m l_{Dj} \left(\left\|x^2-x^1\right\|+\left\|u^2-u^1\right\|\right)\\
 =& \left(\|l_S\|_1+\|l_D\|_1\right)\left(\left\|x^2-x^1\right\|+\left\|u^2-u^1\right\|\right),
\end{align*}
where $l_S=(l_{S1},l_{S2},\ldots, l_{Sm})^\top$ and $l_D=(l_{D1},l_{D2},\ldots, l_{Dm})^\top$, and if
\begin{align*}
\left<G\left(x,u^1\right),u^2-u^1\right> &=\left<S\left(u^1, x\right)-D\left(u^1, x\right),u^2-u^1\right> = \left<S\left(u^1, x\right),u^2-u^1\right> +\left<D\left(u^1, x\right),u^1-u^2\right> \nonumber\\
&=\sum_{j=1}^m S_j\left(u^1, x\right)\left(u_j^2-u_j^1\right) +  \sum_{j=1}^m D_j\left(u^1, x\right)\left(u_j^1-u_j^2\right) \geq 0,
\end{align*}
then
\begin{align*}
\left<G\left(x,u^2\right),u^2-u^1\right>
&= \left<S\left(u^2, x\right)-D\left(u^2, x\right),u^2-u^1\right> =\left<S\left(u^2, x\right),u^2-u^1\right> +\left<D\left(u^2, x\right) ,u^1-u^2\right> \nonumber\\
&=  \sum_{j=1}^m S_j\left(u^2, x\right)\left(u_j^2-u_j^1\right) +  \sum_{j=1}^m D_j\left(u^2, x\right)\left(u_j^1-u_j^2\right) \geq \eta_{SD}\left\|u^2-u^1\right\|^2
\end{align*}
for all $x \in R^n$ and $u^1=(u^1_1,u^1_2,\ldots, u^1_m)^\top, u^2=(u^2_1,u^2_2,\ldots, u^2_m)^\top \in R^m$. Hence $G$ is $\left(\|l_S\|_1+\|l_D\|_1\right)$-Lipschitz and strongly pseudomonotone w.r.t. its second argument, i.e., hypothesis (H$_1$) holds.

Let $u= (u_1,u_2,\ldots, u_m)^\top=(u_j,u_{-j})^\top\in R^m$, $\Omega=\left\{u\in R^m: c_j \leq u_j \leq d_j, j= 1,2,\ldots, m\right\}$ and $\phi(u)=(\phi_1(u_{-1}),\phi_2(u_{-2}),\ldots, \phi_m(u_{-m}))^\top$. Obviously, $\Omega$ is nonempty, closed and convex. In view of hypothesis (A$_2$), we have that (\ref{KlG-2}) holds, and
\begin{equation}\label{Ku-ex}
  K(u)=\left\{u \in R^m: a_j(u_{-j}) \leq u_j \leq b_j(u_{-j}), j = 1,2,\ldots, m\right\}=\phi(u)+\Omega
\end{equation}
with
\begin{align*}
&\left\|\phi\left(u^2\right)-\phi\left(u^1\right)\right\|
\leq \left\|\phi\left(u^2\right)-\phi\left(u^1\right)\right\|_1
=\sum_{j=1}^m \left|\phi_j\left(u_{-j}^2\right)-\phi_j\left(u_{-j}^1\right)\right| \\
\leq& \sum_{j=1}^m l_j \left\|u_{-j}^2 -u_{-j}^1\right\| \leq \left\|u^2 -u^1\right\| \sum_{j=1}^m l_j
= \|l\|_1 \left\|u^2 -u^1\right\|
\end{align*}
for all $u^1=\left(u_j^1,u_{-j}^1\right)^\top, u^2=\left(u_j^2,u_{-j}^2\right)^\top \in R^m$. Consequently, $\phi$ is $\|l\|_1$-Lipschitz. It follows from (\ref{u-PKOmega}) that hypothesis (H$_2$) holds.

For any $y^1,y^2 \in C(I,R^n)$, it follows that
\begin{eqnarray*}
\left\|\psi(y^2)-\psi(y^1) \right\|
  &=& \left\|\sum\limits_{i=1}\limits^k \iota_i y^2(t_i)-\sum\limits_{i=1}\limits^k \iota_i y^1(t_i)\right\|
\leq \sum\limits_{i=1}\limits^k |\iota_i|\left\|y^2(t_i)-y^1(t_i)\right\|\nonumber\\
&\leq& \left(\sum\limits_{i=1}\limits^k |\iota_i|\right) \underset{t\in I}\max\left\|y^2(t)-y^1(t)\right\|
= \left(\sum\limits_{i=1}\limits^k |\iota_i|\right) \underset{t\in I}\max \Big\{e^{\gamma t} e^{-\gamma t} \left\|y^2(t)-y^1(t)\right\|\Big\}\\
&\leq& \left(\sum\limits_{i=1}\limits^k |\iota_i|\right) e^{\gamma T} \left\|y^2-y^1\right\|_\mathcal{B}=l_\psi\left\|y^2-y^1\right\|_\mathcal{B},
\end{eqnarray*}
where $l_\psi=\left(\sum\limits_{i=1}\limits^k |\iota_i|\right) e^{\gamma T}$. Obviously, there exists $\gamma>0$ such that $0<l_\psi <1$. In fact, in view of $0<\sum\limits_{i=1}\limits^k |\iota_i| <1$, taking $\gamma <\frac{1}{T}\ln\frac{1}{\sum\limits_{i=1}\limits^k |\iota_i|}$, one has $0<l_\psi<1$. So hypothesis (H$_5$) holds.

Let $g=\vartheta$ and $f=\chi$. It is clear that hypotheses (H$_3$) and (H$_4$) hold. Note that hypothesis (A$_2$) implies that (\ref{KlG-2}) holds. In light of Theorem \ref{exist-GFDQVI}, one gets that PCP (\ref{GFDQVI-ex}) is unique solvable. \hfill $\Box$
\end{theorem}

\begin{theorem}
If all hypotheses of Theorem \ref{PCP} hold, then PCP (\ref{GFDQVI-ex}) is MLHU stable.

\textbf{Proof.}\hspace{0.2cm}  Applying Theorems \ref{GFDQVI-MLUH} and \ref{PCP}, it follows immediately that PCP (\ref{GFDQVI-ex}) is MLHU stable. \hfill $\Box$
\end{theorem}

\begin{theorem}
Let all hypotheses of Theorem \ref{PCP} hold. If $\varphi \in C(I,R_+)$ is a nondecreasing function, then PCP (\ref{GFDQVI-ex}) is generalized MLHUR stable  w.r.t. $\varphi \mathrm{E}_q$.

\textbf{Proof.}\hspace{0.2cm}  Using Theorems \ref{GFDQVI-GMLUHR} and \ref{PCP}, it follows immediately that PCP (\ref{GFDQVI-ex}) is generalized MLHUR stable w.r.t. $\varphi \mathrm{E}_q$. \hfill $\Box$
\end{theorem}

\section{Conclusions}
\noindent \setcounter{equation}{0}
Throughout this work, we discussed a new system of NFDQVI (\ref{GFDQVI}), which captures the properties of both an FDE with a nonlocal condition and a time-dependent QVI within the same framework. We first showed some properties of solution for PQVI in (\ref{GFDQVI}) under the hypotheses of strong pseudomonotonicity and Lipschitzean. We also showed the unique solvability of NFDQVI (\ref{GFDQVI}) by using the Banach fixed point principle and then obtained some H-U stability results for NFDQVI (\ref{GFDQVI}).  Finally, the abstract results obtained in this work are applicable to two practical problems concerning an MAOP and a PCP.

We would like to mention that the stochastic DVI can be applied to a great deal of real problems arising in many fields such as economy, finance and mechanics in stochastic environments \cite{wang23,zhang2023r,zhang2023s,zhang2023,ZZH,ZZH25}. Thus, it would be important and interesting to investigate some new fractional stochastic DQVI systems under some mild conditions. This is the direction of our future efforts.

\end{document}